\documentclass[preprint,pre,showpacs]{revtex4-1}

\usepackage{graphicx}
\usepackage{epsfig,amssymb,graphicx,color,lineno,enumitem}
\usepackage{setspace}
\usepackage[fleqn]{amsmath}
\usepackage[usenames,dvipsnames]{xcolor}

\begin{document}
\title{Corrected pair correlation functions for environments with obstacles}
\author{Stuart T. Johnston$^{1,2}$}
\email{stuart.johnston@unimelb.edu.au}
\author{Edmund J. Crampin$^{1,2,3}$}
\affiliation{1. Systems Biology Laboratory, School of Mathematics and Statistics, and Department of Biomedical Engineering, University of Melbourne, Parkville, Victoria 3010, Australia,}
\affiliation{2. ARC Centre of Excellence in Convergent Bio-Nano Science and Technology, Melbourne School of Engineering, University of Melbourne, Parkville, Victoria 3010, Australia,}
\affiliation{3. School of Medicine, Faculty of Medicine Dentistry and Health Sciences, University of Melbourne, Parkville, Victoria 3010, Australia.}

\begin{abstract}
Environments with immobile obstacles or void regions that inhibit and alter the motion of individuals within that environment are ubiquitous. Correlation in the location of individuals within such environments arises as a combination of the mechanisms governing individual behavior and the heterogeneous structure of the environment. Measures of spatial structure and correlation have been successfully implemented to elucidate the roles of the mechanisms underpinning the behavior of individuals. In particular, the pair correlation function has been used across biology, ecology and physics to obtain quantitative insight into a variety of processes. However, na\"ively applying standard pair correlation functions in the presence of obstacles may fail to detect correlation, or suggest false correlations, due to a reliance on a distance metric that does not account for obstacles. To overcome this problem, here we present an analytic expression for calculating a corrected pair correlation function for lattice-based domains containing obstacles. We demonstrate that this corrected pair correlation function is necessary for isolating the correlation associated with the behavior of individuals, rather than the structure of the environment. Using simulations that mimic cell migration and proliferation we demonstrate that the corrected pair correlation function recovers the short-range correlation known to be present in this process, independent of the heterogeneous structure of the environment. Further, we show that the analytic calculation of the corrected pair correlation derived here is significantly faster to implement than the corresponding numerical approach.
\end{abstract}


\maketitle

\section{Introduction}

\begin{figure}
\begin{center}
\includegraphics[width=0.8\textwidth]{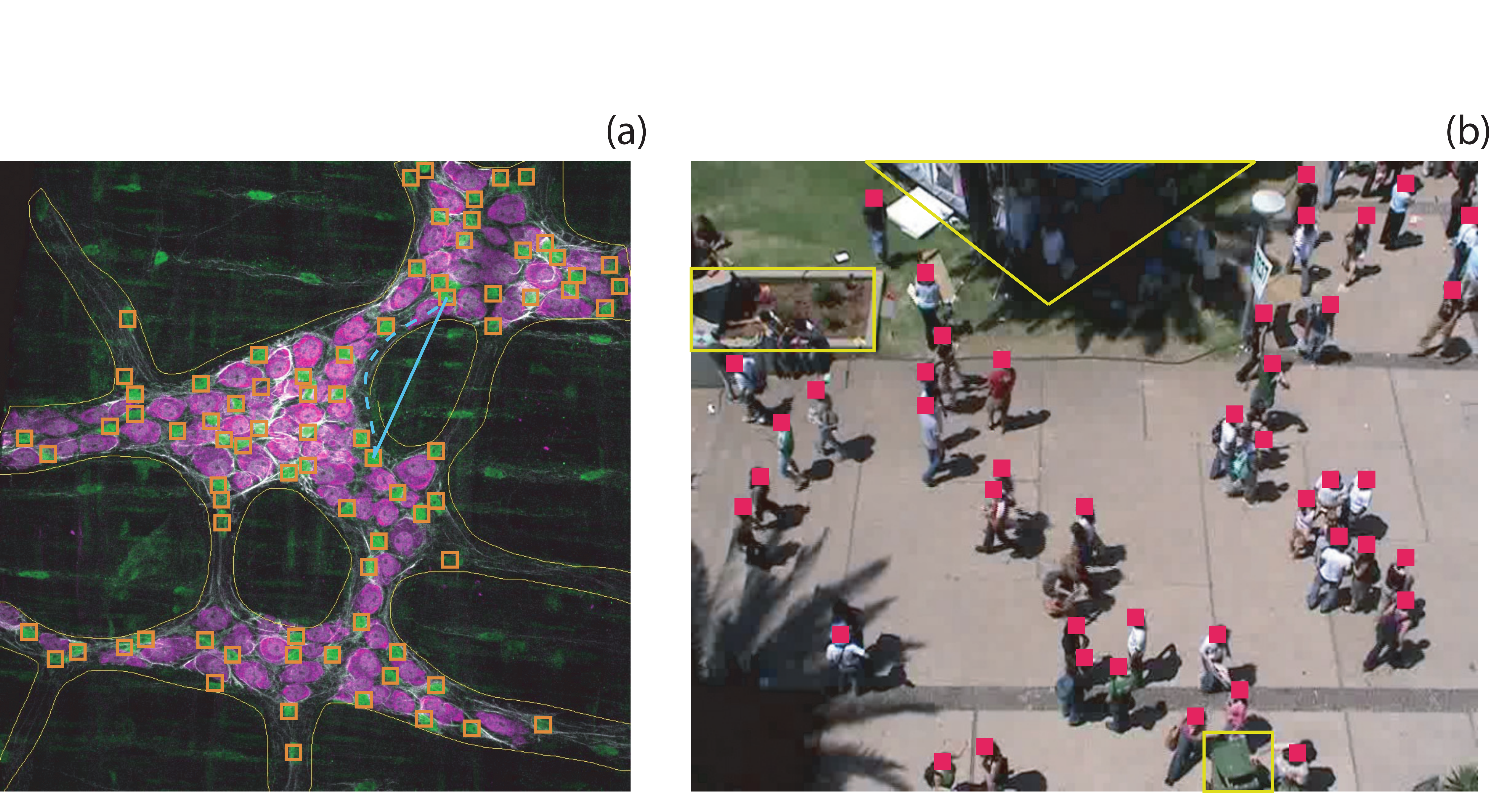}
\caption{(a) Experimental image of the nervous system within the mouse colon, containing neurons (magenta), glial cells (green) and glial processes (white). Glial cells within clusters (known as ganglia) are highlighted with orange squares. Yellow lines indicate inaccessible regions. Path distance and Cartesian distance between cells are highlighted in cyan (dashed and solid, respectively). (b) Experimental image of pedestrian locations (red squares) in the presence of obstacles (yellow lines). Image is obtained from the freely available dataset provided by \cite{lerner2007}.}
\label{F1}
\end{center}
\end{figure}

\noindent Environments that contain obstacles are of interest in a wide variety of fields \cite{dietrich2014,gabella1990,grima2007,helbing2005,hofling2013,lerner2007,moussaid2011,rehder2015,roosen2011,ruhl2005,simpson2017,varas2007}. In biology, it is well known that the motion of macromolecules and proteins in the cytosol is restricted by the densely crowded nature of the interior of cells \cite{hofling2013,roosen2011}. The mesh-like structure of the enteric nervous system, highlighted in Figure \ref{F1}(a), contains clusters of glial cells connected by nerve strands, as well as regions that are inaccessible to the enteric glial cells \cite{gabella1990,ruhl2005}. Hence the location and movement of glial cells is constricted by these inaccessible regions \cite{gabella1990}. In the context of pedestrian dynamics, successful navigation around an obstacle without jamming is a key criteria for the design of safe egress routes \cite{alizadeh2011,dietrich2014,helbing2005,moussaid2011,varas2007}. Similarly, predicting how pedestrians will react to path-blocking obstacles is a key question in computer vision, as developing algorithms for robots to reliably avoid collisions with pedestrians is crucial \cite{rehder2015,ziebart2009}. \\

\noindent Within such environments individuals can undergo self-organization and form highly spatially structured populations, such as the aforementioned clusters of glial cells \cite{gabella1990,ruhl2005} or pedestrian lanes \cite{chraibi2010,moussaid2011}. Quantifying the amount of spatial structure present within an environment provides insight into the mechanisms by which the individuals are governed. Therefore, measures that can be applied to experimental data to obtain estimates of the spatial structure within the data are critical \cite{perry2002}. Various methods for quantifying spatial structure have been proposed previously (for example, see the review by Perry \emph{et al.} \cite{perry2002}, and references therein). Here we focus on the use of pair correlation functions (PCFs), which are a powerful and versatile tool for analysing spatial structure and spatial correlation \cite{agnew2014,bahcall1983,binder2013,binder2015,binny2016,dini2018,holzmann1999,gavagnin2018,johnston2014,johnston2016,law2009}. Pair correlation functions have been successfully employed in astrophysics \cite{bahcall1983}, particle physics \cite{holzmann1999}, ecology \cite{flugge2012,law2009,rajala2018} and cell biology \cite{binder2013,johnston2014}, amongst others. Briefly, the pair correlation for a given distance $m$, $P(m)$, can be defined as
\begin{equation*}
P(m) = \frac{C(m)}{E\left[C(m)\right]},
\end{equation*}
where $C(m)$ is the number of pairs of individuals separated by a distance $m$ observed in the data, and the normalization term, $E\left[C(m)\right]$, is the expected number of individuals separated by distance $m$ if the individuals are located randomly throughout the experimental domain. If there are more individuals separated by a particular distance than expected for randomly located individuals then $P(m) > 1$, and hence there is spatial structure corresponding to aggregation at that distance. Similarly, if fewer individuals are separated by a particular distance than expected then $P(m) < 1$, which suggests that there is spatial structure corresponding to segregation present at distance $m$ \cite{binder2013}. \\

\begin{figure}
\begin{center}
\includegraphics[width=1.0\textwidth]{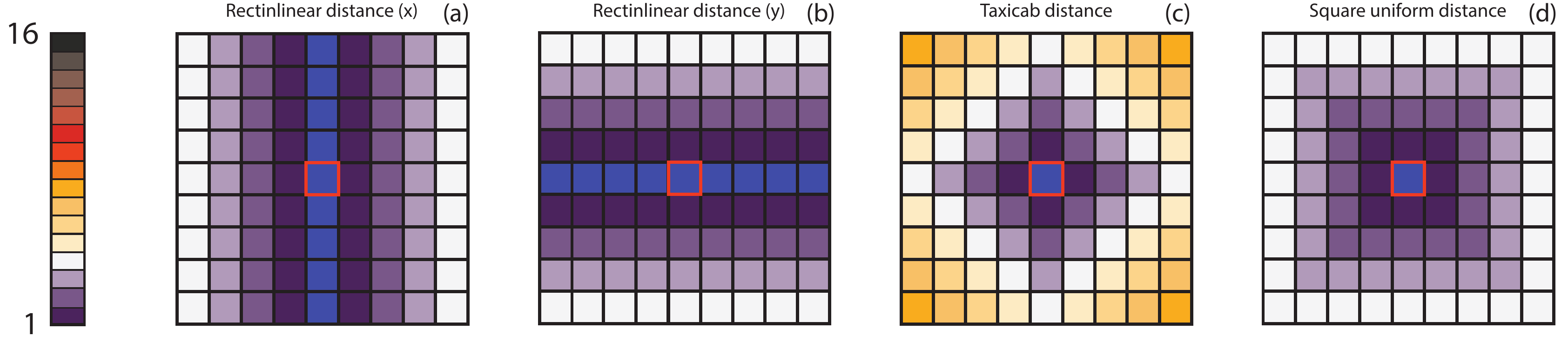}
\caption{Distance between the center lattice site (highlighted in red) and other lattice sites under the (a) rectilinear $x$, (b) rectilinear $y$, (c) taxicab and (d) square uniform distance metrics. Note that blue sites correspond to a distance of zero.}
\label{F2}
\end{center}
\end{figure}

\noindent As it is unlikely that any two pairs of individuals are separated by exactly the same distance, $m$ is typically divided into bins \cite{binder2015}. This can take the form of considering the environment as continuous space and binning the measured distance between pairs of individuals, or by mapping individuals onto a discrete domain such that there is a finite number of possible distances between pairs of individuals \cite{binder2013,binder2015,johnston2014,gavagnin2018}. There has been significant recent focus on PCFs for discrete, or lattice-based, domains \cite{agnew2014,binder2013,binder2015,gavagnin2018,johnston2014,johnston2016}, and deriving analytic expressions for the normalization term under various distance metrics \cite{binder2013,gavagnin2018}. In particular, Binder and Simpson \cite{binder2013} present a normalization term for rectilinear distance in $x$ and $y$, illustrated in Figures \ref{F2}(a)-(b), which corresponds to the distance separating two lattice sites in $x$ and $y$, respectively. More recently, Gavagnin \emph{et al.} \cite{gavagnin2018} derive a normalization term under the taxicab and square uniform distance metrics, illustrated schematically in Figures \ref{F2}(c)-(d), respectively. Under the taxicab distance metric and the square uniform distance metric, the distance between two lattice sites can be thought of as the number of ``jumps'' between the two sites under movement occurring in a von Neumann neighborhood (four nearest-neighbors) and a Moore neighborhood (eight nearest-neighbors), respectively \cite{gavagnin2018}. \\

\noindent However, while these PCFs have proven useful in a range of applications, they are unsuitable for analyzing environments that contain inaccessible regions, due to either ``holes'' in the domain or the presence of obstacles. In these environments, Cartesian distance measures do not adequately describe the distance between two individuals. For example, for the glial cells presented in Figure \ref{F1}(a), certain cells are separated by a path distance that is significantly longer than the Cartesian distance  between the cells. Therefore, na\"ively calibrating a standard PCF to this data may result in a lack of identification of spatial correlation between cells, or may result in spurious correlations being reported. Here we propose an analytic method for calculating a corrected PCF (cPCF) for lattice domains containing obstacles or inaccessible regions. In Section \ref{s:Derivation}, we construct the cPCF in a systematic manner, first considering a single inaccessible site, and subsequently increasing the number of sites within an inaccessible region, as well as increasing the number of inaccessible regions. Through comparison with path-finding algorithms, we show that the derived normalization term is exact. In Section \ref{s:Results}, we demonstrate that this cPCF is required to isolate the correlation associated with the mechanisms governing the behavior of individuals from correlation associated with the structure of the environment. Further, we show that analysis with the exact normalization term is significantly less computationally intensive to perform, compared to using a path-finding algorithm, and we discuss environments where an approximation to the normalisation term can be used effectively. Finally, in Section \ref{s:Discussion}, we discuss our results and suggest potential avenues for future research.

\section{Derivation}
\label{s:Derivation}
\begin{figure}
\begin{center}
\includegraphics[width=1.0\textwidth]{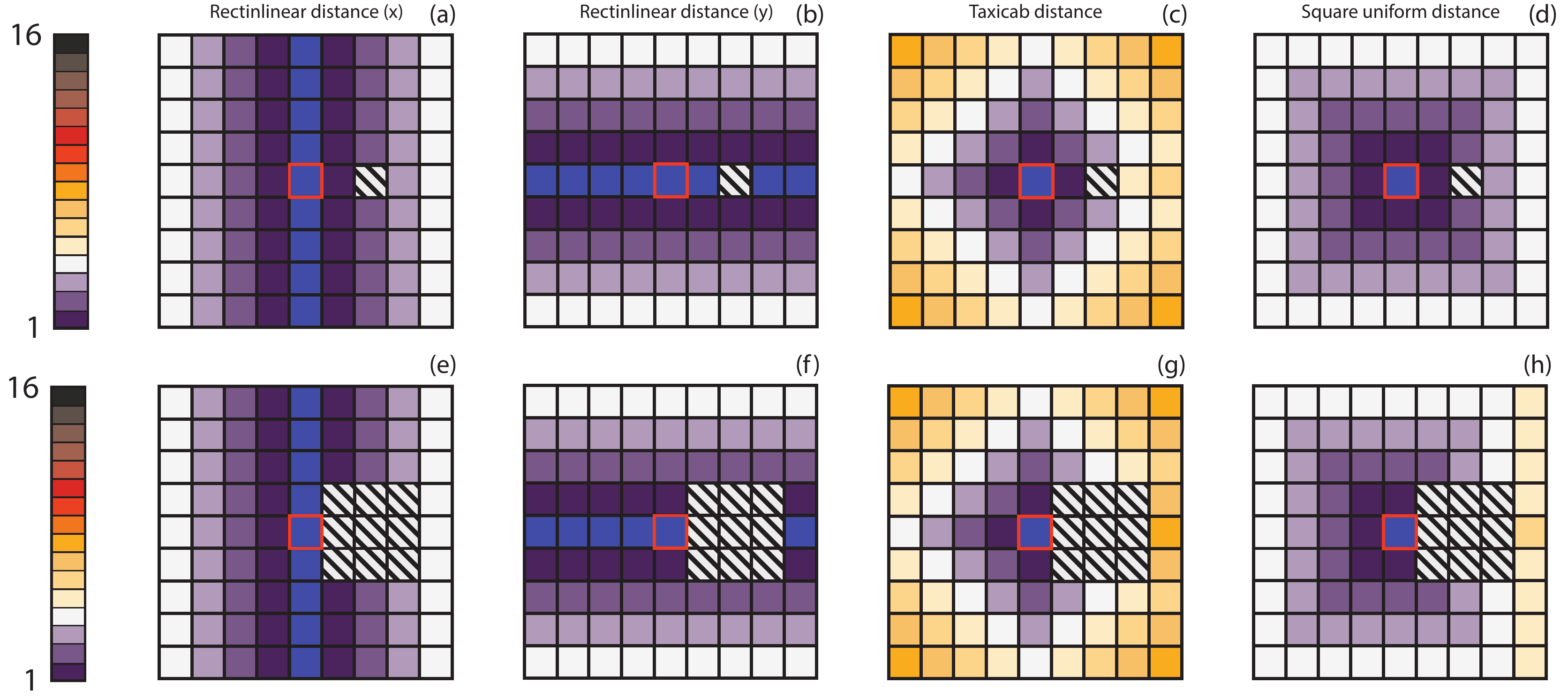}
\caption{Distance between the center lattice site (highlighted in red) and other lattice sites under (a),(e) rectilinear $x$, (b),(f) rectilinear $y$, (c),(g) taxicab and (d),(h) square uniform distance metrics in the presence of (a)-(d) one or (e)-(h) nine inaccessible sites (cross-hatched). Note that blue sites correspond to a distance of zero.}
\label{F3}
\end{center}
\end{figure}

\noindent First we illustrate domains where a different distance metric is required. In Figure \ref{F3}(a)-(d) we introduce a single inaccessible lattice site and examine how the distance metrics change, compared to the domain in Figure \ref{F2}. We note that neither of the rectilinear distances change, and hence all sites remain the same distance from the center site, excluding the inaccessible sites. As such, the counts of pair distances are only reduced by the reduction in the number of lattice sites. The taxicab distance is more significantly impacted by the introduction of the inaccessible site, as to travel from the center of the domain to the rightmost side now requires that the inaccessible site is avoided. Hence the taxicab distance between the center site and sites on the opposite side of the inaccessible site, with respect to the center site, increases by two (Figure \ref{F3}(c)). The square uniform distance is also unaffected by the inaccessible site, as diagonal ``jumps'' count the same as either horizontal or vertical ``jumps''. Therefore the inaccessible site can be avoided by two diagonal ``jumps'', rather than the two horizontal ``jumps'', and the distance does not change (Figure \ref{F3}(d)). If we introduce a larger inaccessible region, as presented in Figures \ref{F3}(e)-(h), we again see that both rectilinear differences are not influenced. As before, the taxicab distance is influenced as the inaccessible region must be avoided to travel between the center site and sites on the right boundary (Figure \ref{F3}(g)). In contrast to the small inaccessible region, the square uniform distance is now affected by the presence of the larger inaccessible region, as the larger region cannot be avoided through diagonal movement (Figure \ref{F3}(h)). We note that any inaccessible region aside from a single site will influence the square uniform distance. As the size of a lattice site typically corresponds to the size of an individual, a distance metric that is not impacted by the presence of obstacles of that size is not appropriate. Furthermore, the majority of models that are implemented on a lattice with obstacles typically only allow movement to one of four nearest-neighbors \cite{alizadeh2011,burstedde2001,ellery2015,ellery2016,ellery2016b,nicolau2007,wedemeier2009}, and hence the taxicab distance metric is implicitly applied. As such, in the remainder of this work, we consider only the taxicab distance metric. \\

\subsection{Standard pair correlation functions}

\noindent To obtain the cPCF for an environment with obstacles under the taxicab distance metric, we first introduce the counts of pair distances for an environment without obstacles. For a domain containing $L_x$ sites in the $x$ direction and $L_y$ sites in the $y$ direction with no-flux boundary conditions, the maximum pair distance  is $L_x+L_y-2$. As experimental images are typically captured such that the influence of boundary effects are minimized, no-flux boundary conditions are perhaps the most relevant boundary conditions \cite{johnston2014}. The alternative choice is periodic boundary conditions; however, we do not expect an individual passing through one boundary to arrive at the opposing boundary. Recently, Gavagnin \emph{et al.} \cite{gavagnin2018} derived the counts of pair distances for a domain without obstacles, $D_{\text{NO}}(m)$, for $m < \text{min}(L_x,L_y)$
\begin{equation*}
D_{\text{NO}}(m) = 2mL_xL_y - (L_x+L_y)m^2 + \frac{m^3-m}{3}.
\end{equation*}
Introducing inaccessible sites into the domain increases the distances between pairs of sites (Figure \ref{F3}), and hence we require an expression for the counts of pair distances for $m \geq \text{min}(L_x,L_y)$. For an arbitrary $L_x$ by $L_y$ domain with no obstacles, the counts of pair distances are (see Appendix A for the derivation)
\begin{equation}
D_{\text{NO}}(m) = \begin{cases} 2mL_xL_y - (L_x+L_y)m^2 + \frac{m^3-m}{3}, \ & 1 \leq m \leq \text{min}(L_x,L_y) \\
				D_{\text{NO}}(\text{min}(L_x,L_y)) - \text{min}(L_x,L_y)^2(m-\text{min}(L_x,L_y)), \ & \text{min}(L_x,L_y) < m < \text{max}(L_x,L_y) \\
				\frac{k(k+1)(k+2)}{3}, \ \text{where} \ k = L_x+L_y-1-m, \ &\text{max}(L_x,L_y) \leq m \leq L_x+L_y-2\end{cases}
				\label{eq:Distances}
\end{equation}

\begin{figure}
\begin{center}
\includegraphics[width=1.0\textwidth]{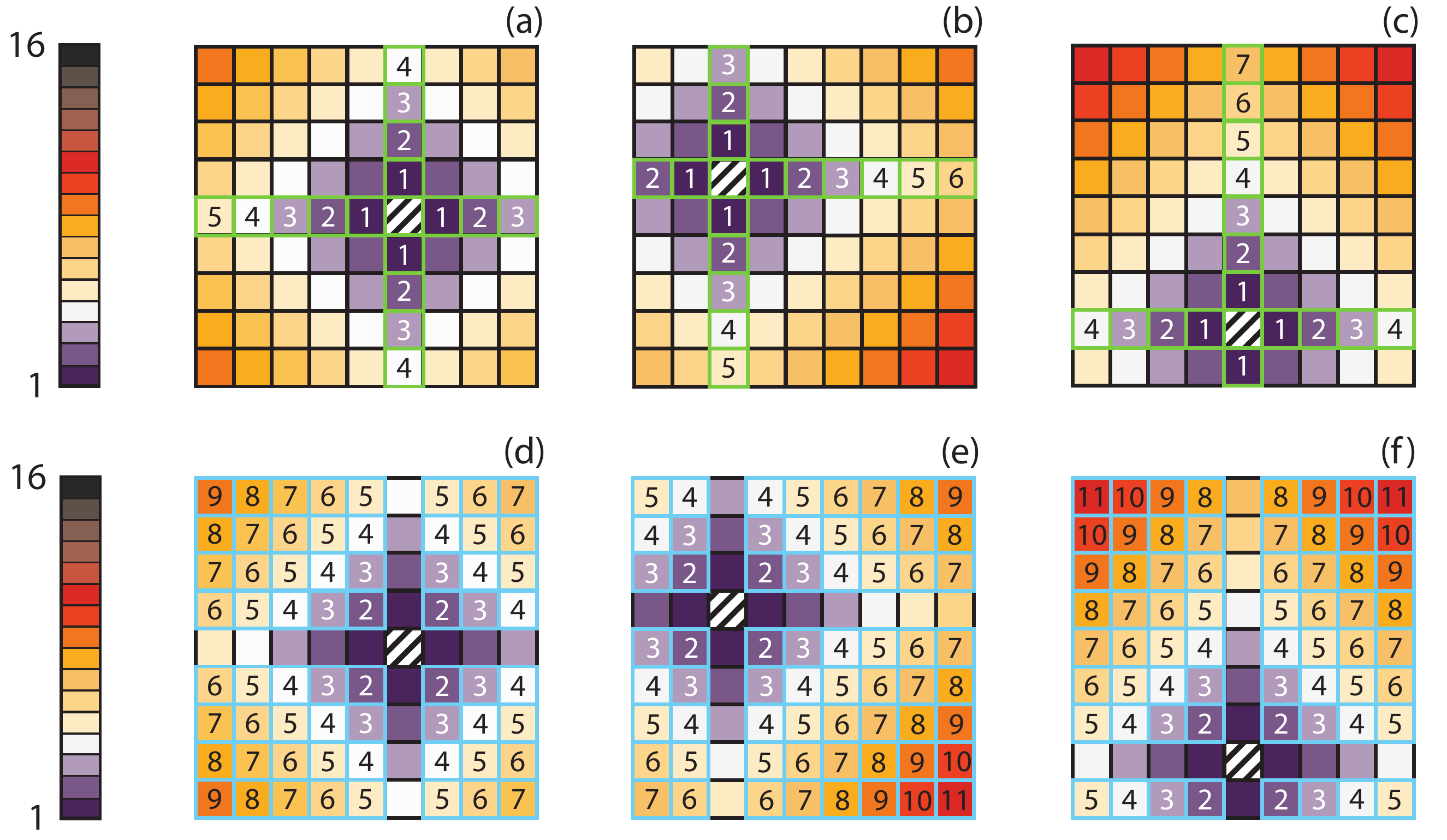}
\caption{Example domains with a single inaccessible site (cross-hatched). The color of individual sites corresponds to the distance between that site and the inaccessible site. Accessible-inaccessible pairs are highlighted in (a)-(c) green or (d)-(f) cyan, depending on whether the accessible site are within the same row or column as the inaccessible site.}
\label{F4}
\end{center}
\end{figure}

\noindent The PCF is calculated by evaluating the counts of pair distances between occupied sites, $C(m)$, and normalizing by the expected number of pair distances obtained from $D(m)$ and the average occupancy of the domain. If there are $z$ occupied sites and $n_a = L_xL_y - n_h$ accessible sites, where $n_h$ is the number of inaccessible sites, then the expected counts of pair distances is \cite{binder2013}
\begin{equation}
E[C(m)] = \frac{z(z-1)}{n_a(n_a-1)}D(m).
\label{eq:PCF}
\end{equation}
Picking two accessible sites at random, $z/n_a$ is the probability that the first selected site is occupied, and $(z-1)/(n_a-1)$ is the probability that the second site is occupied, given that an occupied site has been selected previously. \\

\noindent There are two counts that must be obtained from the data: the counts of pair distances between occupied sites, $C(m)$, and the counts of pair distances between accessible sites, $D(m)$. While $C(m)$ may have to be obtained via a path-finding algorithm, the number of occupied sites is typically small compared to the total number of sites. As such, calculating $C(m)$ will require significantly fewer iterations of the path-finding algorithm, which scales with the square of occupied sites for $C(m)$ or the square of accessible sites for $D(m)$ \cite{gavagnin2018}. Hence, even if calculating $D(m)$ via a path-finding algorithm is prohibitively computationally intensive, $C(m)$ should be able to be calculated rapidly. 

\subsection{Corrected pair correlation functions}
\noindent Next, we focus on obtaining an expression for $D(m)$ for domains containing inaccessible sites, by adjusting the counts of pair distances for a domain with no obstacles, $D_{\text{NO}}$, to account for inaccessible sites. This takes the form of several additional terms:
\begin{itemize}
\item \textit{Accessible-inaccessible pairs}, denoted $A(m)$, which are pairs of sites in the domain that consist of an accessible site and an inaccessible site. As these pairs are counted in $D_{\text{NO}}(m)$, we require that $A(m)$ is accounted for via the removal of these pairs.
\item \textit{Inaccessible-inaccessible pairs}, denoted $I(m)$, which are pairs of inaccessible sites. Again, these sites are counted in $D_{\text{NO}}(m)$ and must be removed to obtain $D(m)$.
\item \textit{Shifted pairs}, denoted $S(m)$, which are pairs of accessible sites where the path distance between the sites is different to the Cartesian distance due to the presence of inaccessible sites. Shifted pairs consist of \textit{lost pairs}, which are pairs of sites in $D_{\text{NO}}(m)$ that are no longer present due to inaccessible sites altering the distance (Figures \ref{F5}(b),(d)), and \textit{gained pairs}, which are pairs of sites that are not in $D_{\text{NO}}(m)$ but are now present due to the introduction of inaccessible sites (Figures \ref{F5}(c),(e)).
\end{itemize}

\subsection{Single inaccessible site}

We will derive an expression for each of the adjustment terms by systematically considering different configurations of inaccessible sites. We first consider a single inaccessible site, such as presented in Figure \ref{F4}, and use the subscript $s$ to denote the special case of a single inaccessible site. The coloring on each site in Figure \ref{F4} highlights the distance between that site and the inaccessible site. For a single inaccessible site, the number of sites with a specific color therefore corresponds to the accessible-inaccessible pairs, $A_s(m)$. In the absence of boundaries, there are 4$m$ pairs for a distance $m$. \\

\noindent Intuitively, the boundaries reduce the number of pairs of sites separated by larger values of $m$. To calculate which of the 4$m$ pairs lie outside the boundary, we introduce eight values, which represent the distance between the inaccessible site and the boundaries and corners of the domain. The distance to the boundary, $b_i$, where $i \in \{L,R,D,U\}$ for the boundary in the left, right, down and up directions, respectively, is defined as
\begin{equation*}
b_L = H_x, \ \ \ \ \  b_R = L_x - H_x + 1, \ \ \ \ \  b_D = H_y,  \ \ \ \ \  b_U = L_y - H_y + 1,
\end{equation*}
where $H_x$ and $H_y$ correspond to the $x$ and $y$ location of the inaccessible site. Similarly, the distance to the corner, $c_{j,k}$ where $j \in \{D,U\}$ and $k \in \{L,R\}$ for the corner in the down-left, down-right, up-left, up-right directions, respectively, is defined as
\begin{equation*}
c_{D,L} = b_D + b_L, \ \ \ \ \  c_{D,R} = b_D + b_R, \ \ \ \ \  c_{U,L} = b_U + b_L,  \ \ \ \ \  c_{U,R} = b_U + b_R.
\end{equation*}
These values allow us to define a function for the number of pairs containing a site that is located outside of the boundaries, $\alpha(m)$, referred to as out-of-domain pairs, and hence
\begin{equation}  
A_s(m) = 4m - \alpha(m).
\label{eq:As}
\end{equation}
The accessible site belonging to an accessible-inaccessible pair can either be located in the same row or column as the inaccessible site (Figures \ref{F4}(a)-(c)), or not located in the same row and not in the same column as the inaccessible site (Figures \ref{F4}(d)-(f)). For sites in either the same row or column as the inaccessible site, there is at most one site at a distance $m$ in each direction (Figures \ref{F4}(a)-(c)). Further, the sites will only be in the domain for distances less than the distance to the boundary. We therefore introduce the function
\begin{equation*}
N_i(m) = \begin{cases} 0, \ &\text{for} \ m < b_i ,\\ 1, \ &\text{for} \ m \geq b_i, \end{cases}
\end{equation*}
which represents the number of out-of-domain pairs of distance $m$ with respect to the boundary $b_i$. For the number of out-of-domain pairs of distance $m$ in diagonal directions, intuitively there will be no pairs if the distance is less than the minimum distance to a boundary of interest. For all distances greater than the minimum distance to a boundary, but less than the distance to the other boundary of interest, we observe that there is the same number of sites in the domain. For example, in Figure \ref{F4}(e), in the down-left direction, we observe only two sites highlighted in cyan for $m = 3,4,5,6$. In comparison, in the down-right direction, where $m = 6$ is less than the minimum distance to a boundary, we observe two, three, four and five sites highlighted in cyan. Hence, the number of out-of-domain pairs increases exactly with distance for distances greater than the minimum distance to a boundary, but less than the distance to the other boundary of interest. For distances greater than the maximum distance to a boundary of interest, we observe that the number of sites highlighted in cyan decreases exactly with distance. Hence the number of out-of-domain pairs increases by two for an increase in distance of one. Finally, for distances greater than the distance to the corner site, there will be no pairs inside the domain, and the number of out-of-domain pairs must be $m-1$. We note that in all cases the sum of number of the pairs outside and inside the domain is $m-1$ for a distance $m$ in a particular diagonal direction. Combining these observations, we introduce the function
\begin{equation*}
M_{j,k}(m) = \begin{cases} 0, \ &\text{for} \ m \leq \text{min}(b_j,b_k), \\ m - \text{min}(b_j,b_k), &\text{for} \ \text{min}(b_j,b_k) < m \leq \text{max}(b_j,b_k), \\ 2m - c_{j,k}, &\text{for} \  \text{max}(b_j,b_k) < m \leq c_{j,k} - 2, \\
m - 1, & \text{for} \ m > c_{j,k} - 2, \end{cases}
\end{equation*}
which represents the number of out-of-domain pairs in each diagonal direction. Note that both $N_i(m)$ and $M_{j,k}(m)$ do not need to be evaluated for pair distances greater than the largest distance between the inaccessible site and the boundary or corner, respectively. The number of out-of-domain pairs is therefore
\begin{equation*}
\alpha(m) = \sum_i N_i(m) + \sum_j \sum_k M_{j,k}(m).
\end{equation*} \\

\begin{figure}
\begin{center}
\includegraphics[width=0.8\textwidth]{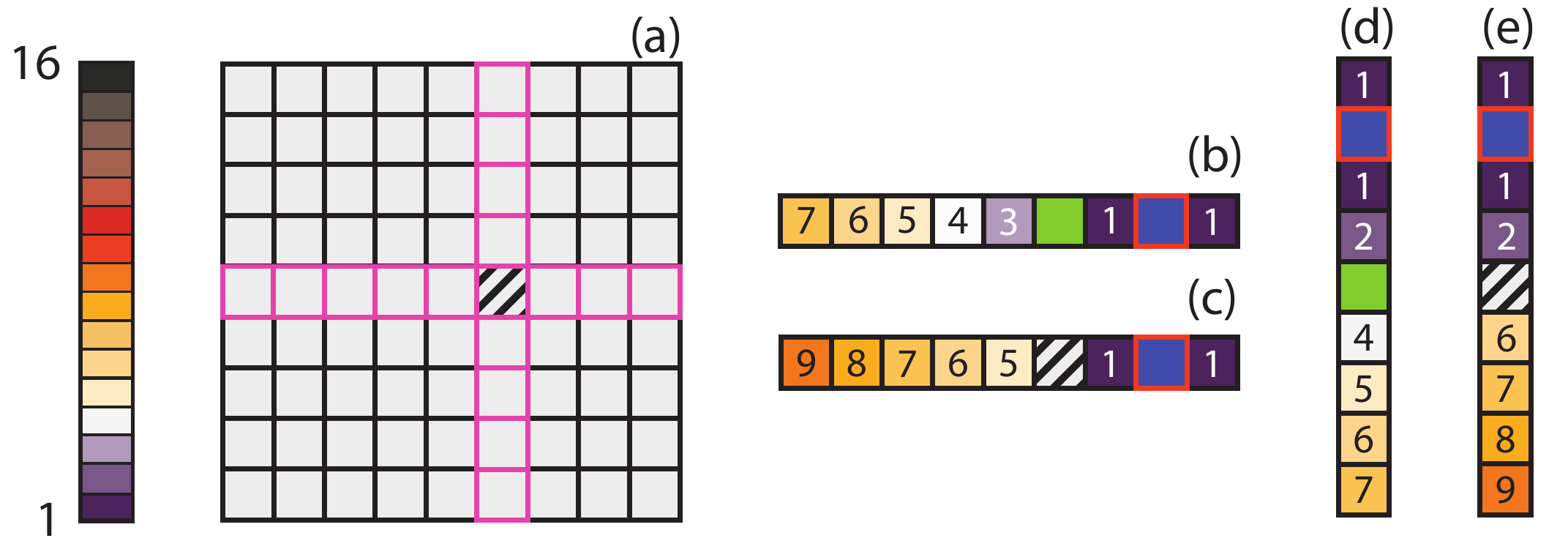}
\caption{(a) Example domain with a single inaccessible site (cross-hatched) and sites contributing to shifted pairs highlighted in pink. Shifted pairs consist of a pair of pink sites, where there is one site on either side of the inaccessible site. (b)-(e) Distance between an example lattice site (highlighted in red) and other lattice sites for (b)-(c) the row of pink sites and (d)-(e) the column of pink sites if the inaccessible site (b),(d) does not have to be avoided (green) or (c),(e) must be avoided.}
\label{F5}
\end{center}
\end{figure}

\noindent For a single inaccessible site there are no inaccessible-inaccessible pairs. Therefore, we next consider the shifted pairs for a single inaccessible site. In this case, pairs of sites that are located on either side of the inaccessible site are shifted pairs, as highlighted in Figure \ref{F5}. The relevant row and column are presented in Figures \ref{F5}(b)-(e), with an example site highlighted to demonstrate the difference between the pair distances in the domain if the inaccessible site is included (Figures \ref{F5}(c) and \ref{F5}(e)) or not (Figures \ref{F5}(b) and \ref{F5}(d)). Intuitively, we observe that if both sites in the pair are on the same side of the inaccessible site then the presence of the inaccessible site is irrelevant. If the sites in the pair are on opposite sides of the inaccessible site, then the inaccessible site increases the pair distance by two, which accounts for a path that avoids the inaccessible site. As we initially consider the counts of pair distances for the domain without inaccessible sites, we can incorporate the shifted pairs by removing the pairs present in Figures \ref{F5}(b) and \ref{F5}(d) and including the pairs present in Figures \ref{F5}(c) and \ref{F5}(e). We note that these are the lost pairs and gained pairs, respectively. \\


\noindent We therefore require an expression for the counts of pair distances within the subdomains in Figures \ref{F5}(b)-(e) for pairs with an inaccessible site on both sides of the inaccessible site. Two observations are useful here; the total number of pairs is conserved, and avoiding the inaccessible site is equivalent to extending the subdomain by two sites, adding two more inaccessible sites, and calculating the pair distances as if the inaccessible site can be passed through. The minimum pair distance for pairs that are located on separate sides of an inaccessible site is one greater than the number of inaccessible sites between them. Hence $L(m)$ is defined on $2 \leq m \leq L_y$ and $2 \leq m \leq L_y$, respectively for the subdomains in Figures \ref{F5}(b) and \ref{F5}(d). The maximum number of pairs separated by a given distance is restricted by the requirement that sites are located on both sides of the inaccessible site, and hence has an upper bound of $d_H = \min(b_L,b_R)-1$ in the horizontal direction and $d_V = \min(b_U,b_D)-1$ in the vertical direction. These values represent the minimum of the number of sites either side of the inaccessible site. Further, there is only one possible pair of sites separated by a distance of 2 and $L_x-1$ (or $L_y-1$), two possible pairs of sites for a distance of 3 and $L_x-2$ (or $L_y-2$), and so forth. Combining this with the aforementioned upper bound we obtain an expression for $L(m)$ for a single inaccessible site,
\begin{equation}
L_s(m) = \text{min}\left(-\left|m-\frac{L_y+1}{2}\right|+\frac{Ly-1}{2},d_V\right) + \text{min}\left(-\left|m-\frac{L_x+1}{2}\right|+\frac{Lx-1}{2},d_H\right).
\label{eq:Ls}
\end{equation}
As noted previously, the number of shifted pairs is conserved. As the pair distance increases by two in the presence of a single inaccessible site, here
\begin{equation}
G_s(m) = \text{min}\left(-\left|m-\frac{L_y+5}{2}\right|+\frac{Ly-1}{2},d_V\right) + \text{min}\left(-\left|m-\frac{L_x+5}{2}\right|+\frac{Lx-1}{2},d_H\right).
\label{eq:Gs}
\end{equation}
This corresponds to an increase of two in both the number of inaccessible sites and the subdomain length. We have now considered all the adjustment terms for a single inaccessible site and therefore by combining \eqref{eq:Distances}, \eqref{eq:As} and \eqref{eq:Ls}-\eqref{eq:Gs}, we obtain the expression for the count of pair distances for a domain with a single inaccessible site, noting $I(m) = 0$ for a single inaccessible site,
\begin{equation}
D_s(m) = D_{\text{NO}}(m) - A_s(m) - L_s(m) + G_s(m).
\end{equation}

\subsection{Clusters of inaccessible sites}

\begin{figure}
\begin{center}
\includegraphics[width=0.8\textwidth]{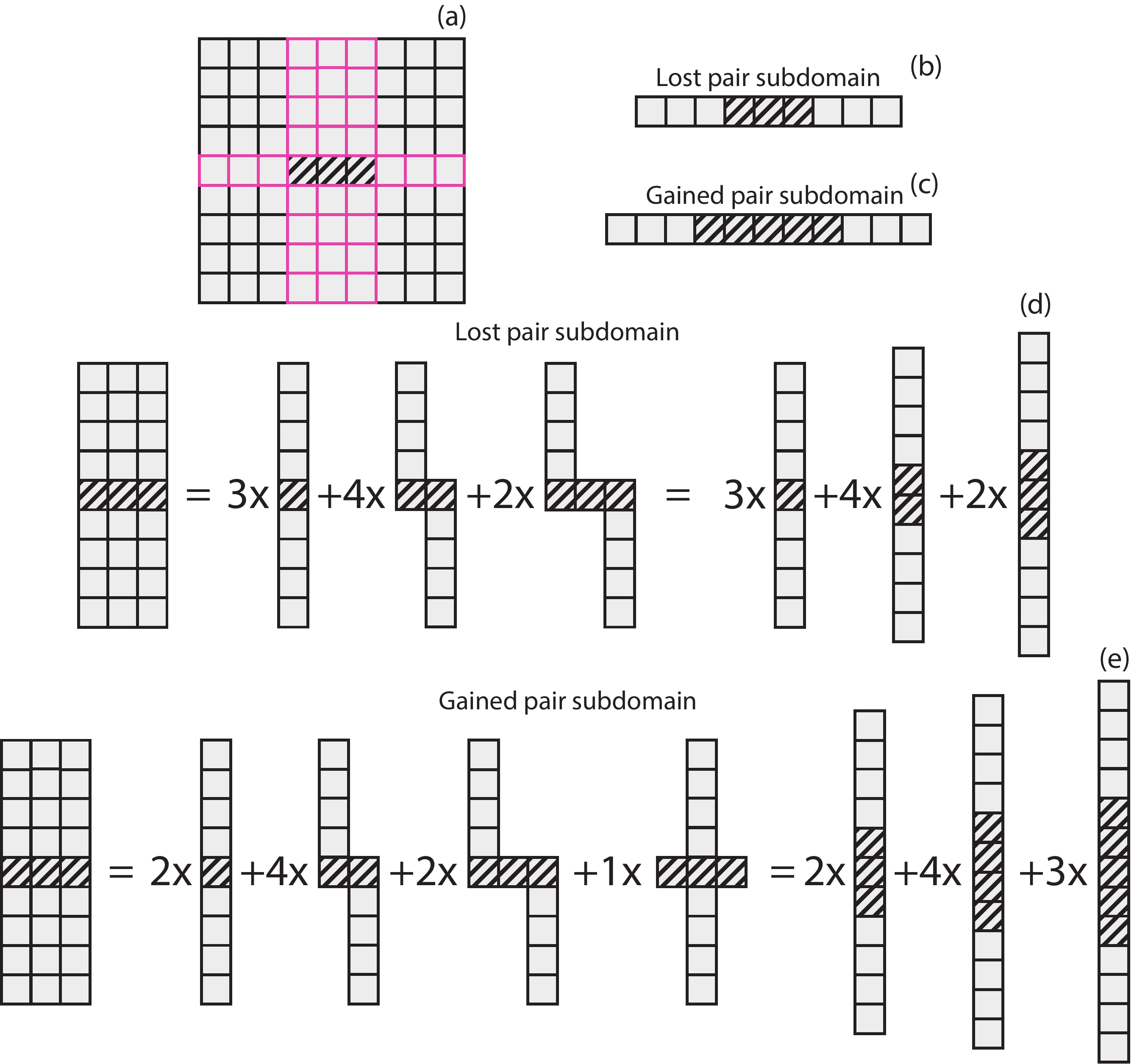}
\caption{(a) Example domain with a cluster of inaccessible sites (cross-hatched) and sites contributing to shifted pairs highlighted in pink. (b) Subdomain that contains all lost pairs for the row of pink sites in (a). (c) Subdomain that contains all gained pairs for the row of pink sites in (a). (d) Subdomain, and associated transformation to multiple subdomains, which contain all lost pairs for the three columns of pink sites in (a). (e) Subdomain, and associated transformation to multiple subdomains, which contain all gained pairs for the three columns of pink sites in (a).}
\label{F6}
\end{center}
\end{figure}

\noindent We now consider a domain with several inaccessible sites that form a $C_x$ by 1 cluster of inaccessible sites, as presented in Figure 6. Introducing the two additional sites in this example results in an increase in the number of accessible-inaccessible pairs. The total number of accessible-inaccessible pairs, $A_s(m)$, is given by
\begin{equation}
A(m) = \sum_H A_s(m),
\label{eq:A}
\end{equation}
where $H$ is the set of inaccessible sites in the domain and $A_s(m)$ is defined in \eqref{eq:As}. Including the additional inaccessible sites introduces inaccessible-inaccessible pairs. The calculation of the number of inaccessible-inaccessible pairs, $I(m)$, is straightforward, as the distance between such pairs is simply the taxicab distance between the two relevant sites:
\begin{equation}
I(m) = \sum_{i = 1}^{n_h}\sum_{j = i+1}^{n_h} \textbf{1}_m\big(d_{\text{taxicab}}(\textbf{h}_i - \textbf{h}_j)\big),
\label{eq:I}
\end{equation}
where $\textbf{1}_m(x)$ is the indicator function, one when $x = m$ and zero otherwise, $n_h$ is the number of inaccessible sites, $d_{\text{taxicab}}(\textbf{h}_i - \textbf{h}_j)$ is the taxicab distance between two sites located at sites $\textbf{h}_i = (H_x^i,H_y^i)$ and $\textbf{h}_j = (H_x^j,H_y^j)$. We observe that $I(m) = 0$ for $n = 1$, as discussed previously. \\

\noindent To calculate lost and gained pairs for a cluster of inaccessible sites, we consider the example domain in Figure \ref{F6}. The lost pairs and gained pairs associated with the horizontal direction are relatively straightforward, and can be calculated from the subdomains presented in Figures \ref{F6}(b)-(c). We introduce the general counts of pair distances for a one-dimensional subdomain $K(m,n,X,d)$, where $X$ is the total number of sites in the subdomain, $n$ is the number of inaccessible sites in the subdomain and $d$ is the relevant $d_V$ or $d_H$ value, as defined previously. The function is
\begin{equation}
K(m,n,X,d) = \min\left(-\left| m - \frac{X + n}{2}\right| + \frac{X - n}{2},d\right).
\end{equation}
We note that this generalizes Equations \eqref{eq:Ls}-\eqref{eq:Gs} by allowing for an arbitrary number of inaccessible sites. By using this generalized definition, we can consider transformations of the cluster of inaccessible sites to multiple one-dimensional subdomains with varying $X$ and $n$, and calculating the $K$ value for each. For example, for the horizontal subdomains in Figure \ref{F6}(b), the lost pair subdomain uses the $X$, $n$ and $d$ values obtained from the original domain, whereas the gained pair subdomain in Figure \ref{F6}(c) uses $X+2$, $n+2$ and $d$. This is consistent with the previous observation of a longer pair distance corresponding to the path around the inaccessible sites. We reiterate that this implies that the sites around the inaccessible sites are accessible. The vertical subdomains in Figure \ref{F6}(a) are more complicated due to the increase in the number of columns where pairs of sites can be located such that the distance between them is influenced by the inaccessible sites. These columns are highlighted in pink in Figure \ref{F6}(a). As the accessible sites must be located on either side of the inaccessible sites, there are a total of nine combinations of columns. Hence, the transformation of the vertical subdomain around the cluster of inaccessible sites results in nine one-dimensional subdomains. For the lost pairs, there can be one, two or three inaccessible sites between the pair of accessible sites. The number of each of these possibilities is highlighted in Figure \ref{F6}(d). Note that an increase in the number of inaccessible sites corresponds to an increase in both $X$ and $n$. For the gained pairs, there can again be one, two or three inaccessible sites between the pair of accessible sites. However, the magnitude of the increase in pair distance due to the inaccessible sites depends on which columns the pair belongs to. If either one of the accessible sites is in an outermost column, then the increase in pair distance is two. However, if both accessible sites are in the middle column, then the increase in pair distance is four, as the path from one site to the other must avoid all of the inaccessible sites. In general, the increase in distance is two times the minimum distance between the columns (rows) that the accessible sites belong to, and the columns (rows) on either side of the inaccessible sites. We illustrate the transformation of the vertical subdomain into one-dimensional subdomains for gained pairs in Figure \ref{F6}(e), as well as the relevant number of each subdomain. In general, an $C_x$ by 1 cluster of inaccessible sites results in
\begin{align}
L(m) &= K(m,C_x,L_x,d_H) + C_xK(m,1,L_y,d_V) \nonumber \\ & + \sum_{i = 2}^{C_x} 2(C_x-i+1)K(m,i,L_y+i-1,d_V), \label{eq:Lm1}
\end{align}
and
\begin{align}
G(m) &= K(m,C_x+2,L_x+2,d_H) + C_xK(m,C_x+2,L_y+C_x+1,d_V) \nonumber \\ & + \sum_{i = 1}^{C_x-1} 2iK(m,i+2,L_y+i+1,d_V). \label{eq:Gm1}
\end{align}
The first term in both equations corresponds to the horizontal subdomain, and the second and third terms correspond to the $(C_x)^2$ vertical subdomains. A similar function for the lost and gained pairs would describe a 1 by $C_y$ cluster of inaccessible sites, as this is simply a rotation of the $C_x$ by 1 cluster of inaccessible sites.\\

\begin{figure}
\begin{center}
\includegraphics[width=1.0\textwidth]{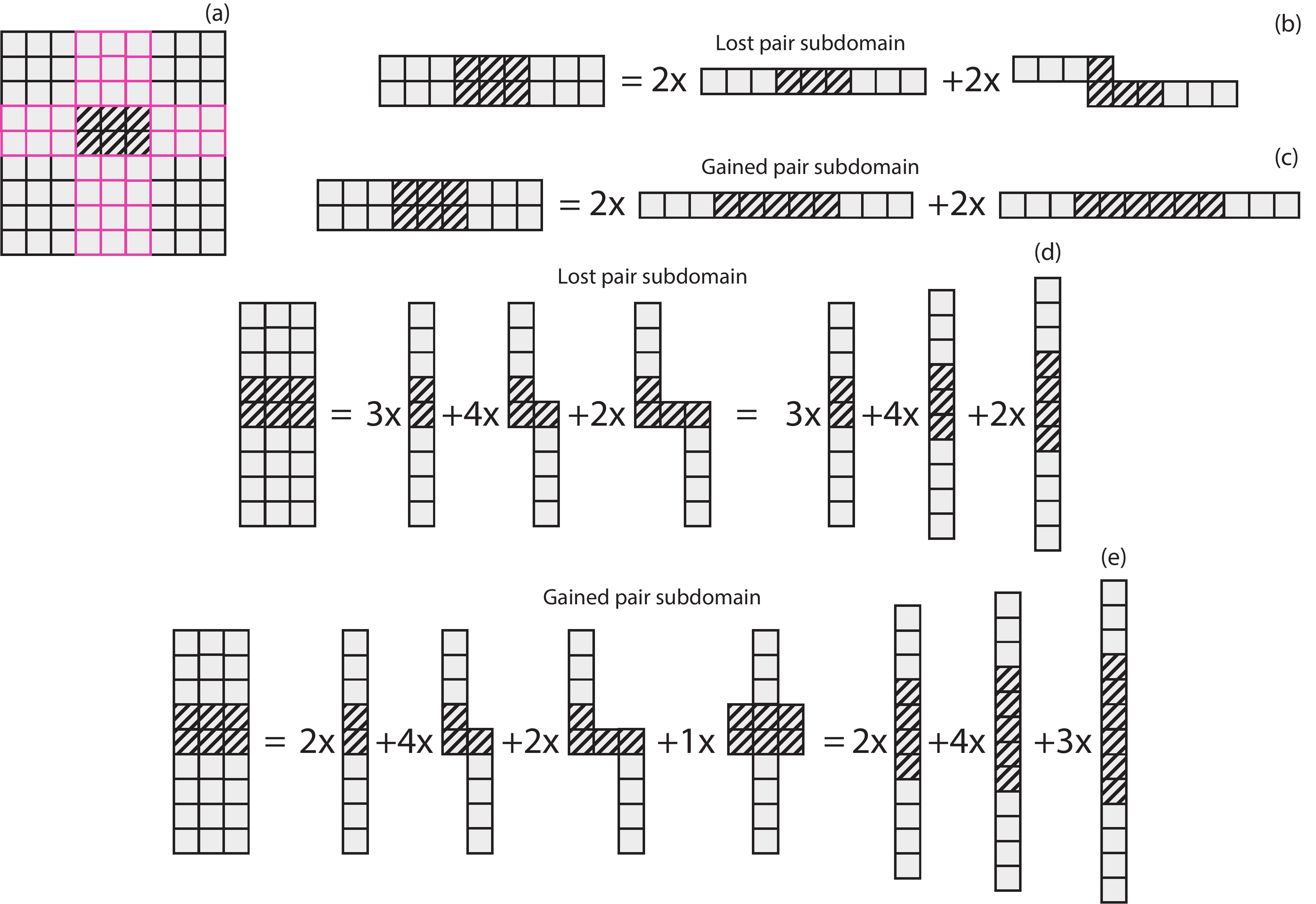}
\caption{(a) Example domain with a cluster of inaccessible sites (cross-hatched) and sites contributing to shifted pairs highlighted in pink. (b)-(c) Subdomain, and associated transformation to multiple subdomains, which contain all lost and gained pairs, respectively, for the two rows of pink sites in (a). (d)-(e) Subdomain, and associated transformation to multiple subdomains, which contain all lost and gained pairs, respectively, for the three columns of pink sites in (a).}
\label{F7}
\end{center}
\end{figure}

\noindent It is now relatively straightforward to extend the lost pairs and gained pairs functions, \eqref{eq:Lm1}-\eqref{eq:Gm1}, to apply to an $C_x$ by $C_y$ cluster of inaccessible sites, as presented in Figure \ref{F7}. We note that the vertical subdomain transformation is similar to the previous $C_x$ by 1 cluster, albeit with an increase in $n$. Now that the cluster has $C_y > 1$, the horizontal subdomain transformation also results in multiple one-dimensional subdomains. We note that the transformation is the same as for the vertical subdomain except for rotation and hence it is straightforward to obtain the lost pairs function for an arbitrary $C_x$ by $C_y$ cluster of inaccessible sites,
\begin{equation}
L_c(m,C_x,C_y) = L_c^h(m,C_x,C_y,L_y) + L_c^v(m,C_x,C_y,L_y), \label{eq:Lc}
\end{equation}
where 
\begin{equation*}
L_c^h(m,C_x,C_y,L_y) = C_yK(m,C_x,L_x,d_H) + \sum_{i = 2}^{C_y} 2(C_y-i+1)K(m,i+C_y-1,L_x+i-1,d_H),
\end{equation*}
is the horizontal contribution to the lost pair distances and
\begin{equation*}
L_c^v(m,C_x,C_y,L_y) = C_xK(m,C_y,L_y,d_V) + \sum_{i = 2}^{C_x} 2(C_x-i+1)K(m,i+C_x-1,L_y+i-1,d_V),
\end{equation*}
is the vertical contribution to the lost pair distances. Similarly, the gained pairs function for an arbitrary $C_x$ by $C_y$ cluster of inaccessible sites can be obtained, and is given by
\begin{equation}
G_c(m,C_x,C_y) = G_c^h(m,C_x,C_y,L_x) + G_c^v(m,C_x,C_y,L_x), \label{eq:Gc}
\end{equation}
where
\begin{equation*}
G_c^h(m,C_x,C_y,L_x) = C_yK(m,C_x+C_y+1,L_x+C_y+1,d_H) + \sum_{i = 1}^{C_y-1} 2iK(m,i+C_x+1,L_x+i+C_x,d_H),
\end{equation*}
is the horizontal contribution to the gained pair distances and
\begin{equation*}
G_c^v(m,C_x,C_y,L_x) = C_xK(m,C_y+C_x+1,L_y+C_x+1,d_V) + \sum_{i = 1}^{C_x-1} 2iK(m,i+C_y+1,L_y+i+C_y,d_V),
\end{equation*}
is the vertical contribution to the gained pair distances. \\

\noindent Combining \eqref{eq:A}-\eqref{eq:I} and \eqref{eq:Lc}-\eqref{eq:Gc} we obtain the number of pair distances for a cluster of inaccessible sites
\begin{equation}
D_c(m) = P_{NO}(m) - A(m) + I(m) - L_c(m) + G_c(m).
\label{eq:Dc}
\end{equation}
Note that the $I(m)$ term is positive as the inaccessible-inaccessible pairs are counted twice and removed in $A(m)$.

\begin{figure}
\begin{center}
\includegraphics[width=0.8\textwidth]{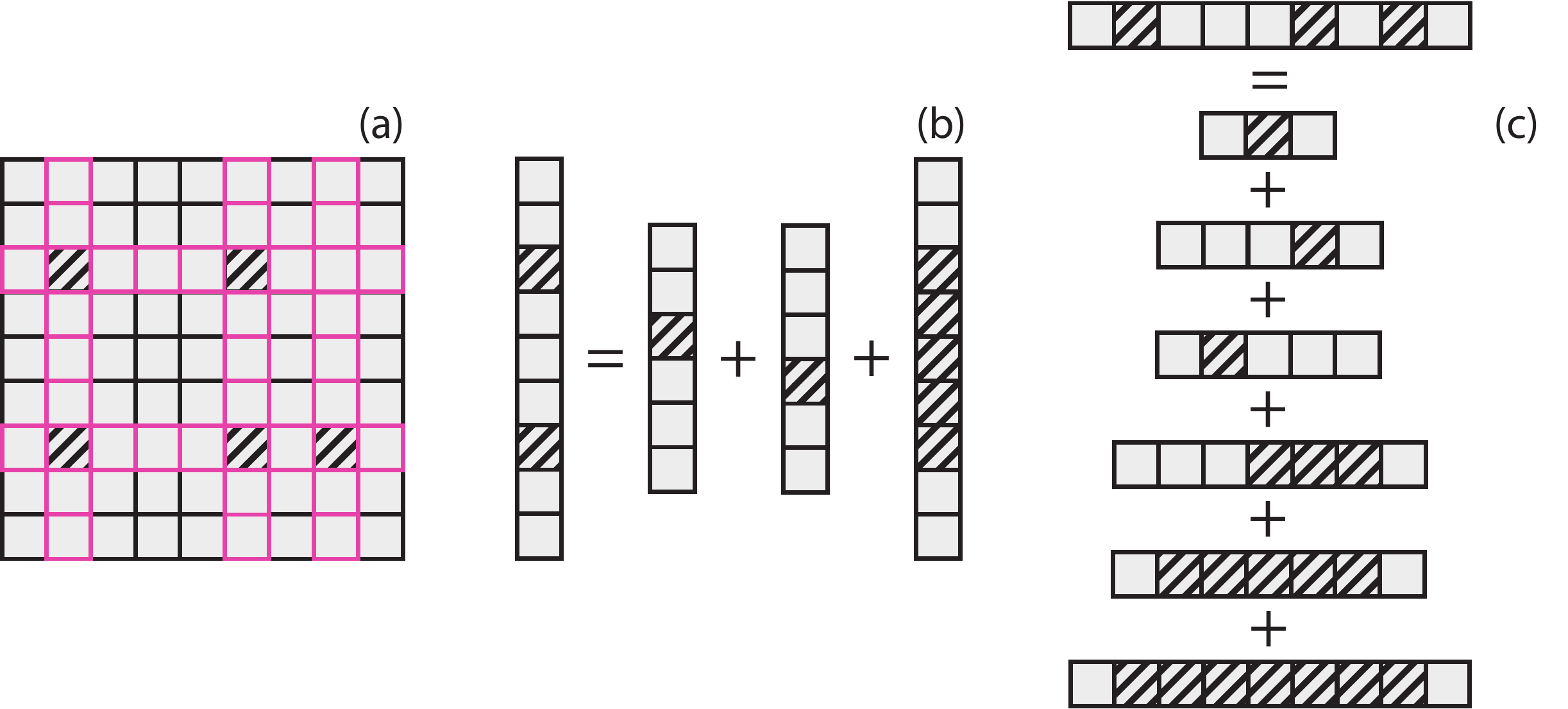}
\caption{(a) Example domain with multiple inaccessible sites (cross-hatched) and sites contributing to shifted pairs highlighted in pink. (b) Subdomain, and associated transformation to multiple subdomains, for the left-most column of pink sites. (c) Subdomain, and associated transformation to multiple subdomains, for the bottom-most row of pink sites. Note that each transformed subdomain has both lost and gained pairs as described previously.}
\label{F8}
\end{center}
\end{figure}

\subsection*{Multiple clusters of inaccessible sites}

\noindent Thus far we have considered only a single cluster of inaccessible sites. We now consider the generalization to multiple clusters of inaccessible sites, as illustrated in Figure \ref{F8}(a). In this example, we consider a domain with 1 by 1 clusters of inaccessible sites. However, we note that the approach generalizes to rectangular clusters of any size provided that the clusters are arranged such that entirely accessible rows and columns exist on either side of all clusters. We also note that the previous definitions of the accessible-inaccessible pairs and inaccessible-inaccessible pairs, \eqref{eq:A}-\eqref{eq:I}, are valid here. Similar to the approach to obtain the lost and gained pairs for a $C_x$ by $C_y$ cluster, a subdomain is isolated and transformed into $B^i(B^i+1)/2$ one-dimensional subdomains, where $B^i$ is the number of clusters within a subdomain. Consider the column subdomain presented in Figure \ref{F8}(b). The transformation isolates all possible combinations of inaccessible sites, and the lost pairs and gained pairs functions are calculated for each one-dimensional subdomain. The transformation of the subdomain in Figure \ref{F8}(b) results in three subdomains. The first two subdomains contain a single inaccessible site and $X$ total sites, where $X$ corresponds to the maximum number of sites from the original subdomain that are contiguous and still contain only that single inaccessible site. The third subdomain contains both inaccessible sites, and renders all of the intervening sites inaccessible. Here $X$ corresponds to the length of the original subdomain. As we are able to calculate the lost and gained pair functions on a subdomain with a single cluster of inaccessible sites, transforming the subdomain with multiple separate inaccessible sites provides a straightforward method for this calculation. We note that the single inaccessible site subdomains provide pairs that contain sites in the upper-most (or bottom-most) region of the original subdomain as well as in the middle of the original subdomain. The subdomain with additional inaccessible sites provides pairs that are located in both the upper-most and bottom-most region of the original subdomain. In Figure \ref{F8}(c), we consider the row with three single inaccessible sites, and the associated transformation. As before, we consider all possible combinations of inaccessible sites: three subdomains with a single inaccessible site, two subdomains with inaccessible regions bounded by two inaccessible sites, and finally a subdomain containing all three inaccessible sites, and the sites between them treated as inaccessible. \\

\noindent We now generalize this approach to $B$ distinct clusters within a row (column) of clusters of inaccessible sites, where $\textbf{s}_{\text{h}}$ ($\textbf{s}_{\text{v}}$) is the set of co-ordinates of the left-most (upper-most) site in inaccessible clusters and $\textbf{f}_{\text{h}}$ ($\textbf{f}_{\text{v}}$) is the set of co-ordinates of the right-most (bottom-most) site in inaccessible clusters.
\begin{equation}
L(m) = \sum_{i=1}^{B} \sum_{j=1}^{B-i+1} L_c^h(m,\textbf{f}_{\text{h}}^{i+j-1}-\textbf{s}_{\text{h}}^{j}-1,C_y,\textbf{s}_{\text{h}}^{i+j} - \textbf{f}_{\text{h}}^{j-1} - 1).
\end{equation}
The first summation represents the number of clusters of inaccessible sites in the domain, and the second summation represents the number of combinations of $i$ neighboring clusters. The $\textbf{f}_{\text{h}}^{i+j-1}-\textbf{s}_{\text{h}}^{j}-1$ value corresponds to the number of sites in the cluster in the horizontal direction, and the $\textbf{s}_{\text{h}}^{i+j} - \textbf{f}_{\text{h}}^{j-1} - 1$ value corresponds to the length of the subdomain. Repeating this processes over all distinct rows and columns of clusters of inaccessible sites, we obtain
\begin{align}
L(m) = & \sum_{i=1}^{n_{\text{rows}}}\sum_{j=1}^{B^i} \sum_{k=1}^{B^i-j+1} L_c^h(m,\textbf{f}_{\text{h},i}^{j+k-1}-\textbf{s}_{\text{h},i}^{k}-1,C_y,\textbf{s}_{\text{h},i}^{j+k} - \textbf{f}_{\text{h},i}^{k-1} - 1) \  + \nonumber \\ & \sum_{i=1}^{n_{\text{columns}}}\sum_{j=1}^{B^i} \sum_{k=1}^{B^i-j+1} L_c^v(m,C_x,\textbf{f}_{\text{v},i}^{j+k-1}-\textbf{s}_{\text{v},i}^{k}-1,\textbf{s}_{\text{v},i}^{j+k} - \textbf{f}_{\text{v},i}^{k-1} - 1),
\label{eq:L}
\end{align}
for the lost pair distances and, following similar arguments, 
\begin{align}
G(m) = & \sum_{i=1}^{n_{\text{rows}}}\sum_{j=1}^{B^i} \sum_{k=1}^{B^i-j+1} G_c^h(m,\textbf{f}_{\text{h},i}^{j+k-1}-\textbf{s}_{\text{h},i}^{k}-1,C_y,\textbf{s}_{\text{h},i}^{j+k} - \textbf{f}_{\text{h},i}^{k-1} - 1) \  + \nonumber \\ & \sum_{i=1}^{n_{\text{columns}}}\sum_{j=1}^{B^i} \sum_{k=1}^{B^i-j+1} G_c^v(m,C_x,\textbf{f}_{\text{v},i}^{j+k-1}-\textbf{s}_{\text{v},i}^{k}-1,\textbf{s}_{\text{v},i}^{j+k} - \textbf{f}_{\text{v},i}^{k-1} - 1),
\label{eq:G}
\end{align}
for the gained pair distances, where $n_{\text{rows}}$ ($n_{\text{columns}}$) is the number of rows (columns) that contain distinct clusters of inaccessible sites. A cluster that has $C_x > 1$ would only contribute once to $n_{\text{columns}}$ rather than $C_x$ times, and we note that the domain in Figure \ref{F8}(a) has $n_{\text{columns}} = 3$ and $n_{\text{rows}} = 2$. \\

\noindent Therefore, the expression for the counts of pair distances for a domain with obstacles is
\begin{equation}
D(m) = D_{\text{NO}}(m) - A(m) + I(m) - L(m) + G(m),
\label{eq:NewDistances}
\end{equation}
where $D_{\text{NO}}(m)$, $A(m)$, $I(m)$, $L(m)$ and $G(m)$ are defined in \eqref{eq:Distances}, \eqref{eq:A}-\eqref{eq:I} and \eqref{eq:L}-\eqref{eq:G}, respectively. Again, we note that this expression is exact provided that the obstacles are arranged such that each obstacle has an entirely vacant row or column on all sides of the obstacle.
\section{Results}
\label{s:Results}

\begin{figure}
\begin{center}
\includegraphics[width=1.0\textwidth]{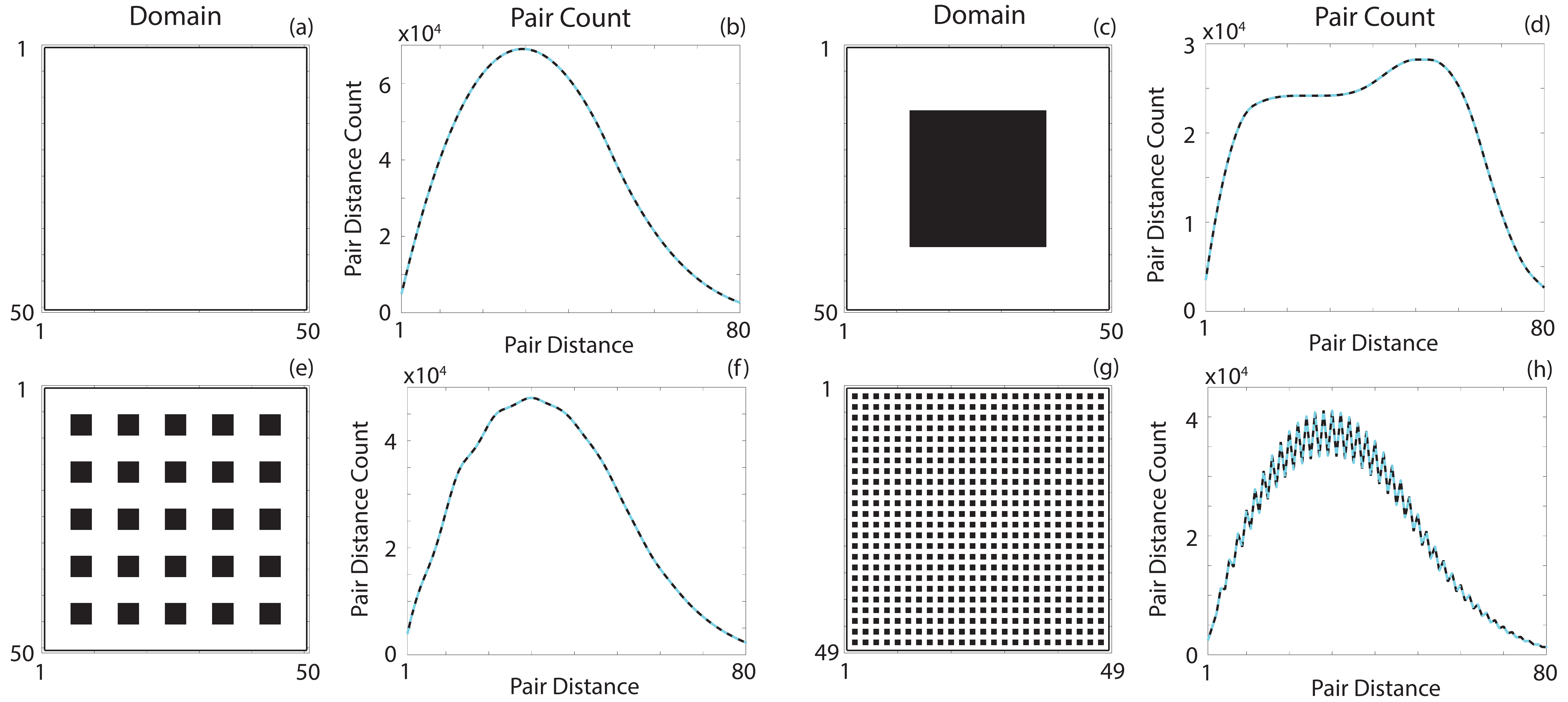}
\caption{(a),(c),(e),(g) Domains containing various configurations of inaccessible sites (black) and (b),(d),(f),(h) the corresponding count of pair distances obtained from the analytic expression (cyan) or numerical approach (black, dashed).}
\label{F9}
\end{center}
\end{figure}

\begin{figure}
\begin{center}
\includegraphics[width=1.0\textwidth]{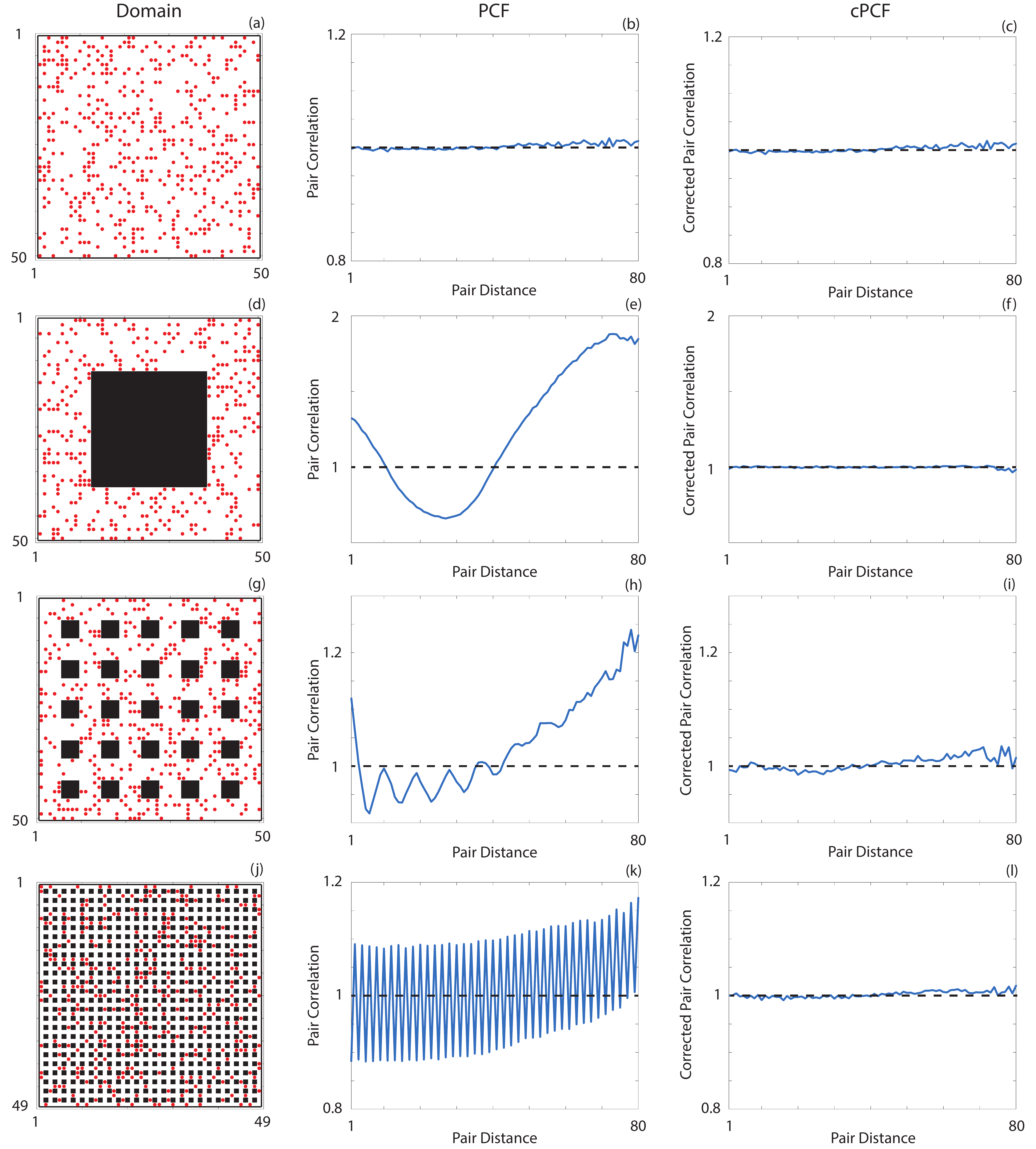}
\caption{(a),(d),(g),(j) Domains containing various configurations of inaccessible sites (black) with agents (red) randomly placed on accessible sites with the corresponding (b),(e),(h),(k) standard PCF, that is, the pair correlation calculated ignoring inaccessible sites, and (c),(f),(i),(l) cPCF. The dashed black line corresponds to no correlation. All PCFs are the average of 100 identically-prepared domains.}
\label{F10}
\end{center}
\end{figure}

\noindent We first verify that the analytic expression \eqref{eq:NewDistances} exactly calculates the counts of pair distances for a domain with obstacles. In Figures \ref{F9}(a),(c),(e),(g) we present four different domains with inaccessible sites highlighted in black. For each domain we calculate the counts of pair distances numerically using Matlab's $\verb"graphshortestpath"$, which calculates the shortest distance between any two points on a graph, given the adjacency matrix of the graph. This is the approach suggested by Gavagnin \emph{et al.} \cite{gavagnin2018} for calculating general PCFs, who note that computational cost is $\mathcal{O}(L^3(2L-2))$ for a square domain with $L$ sites in each direction, which can become computationally infeasible even for modest $L$ values . We then evaluate the analytic expression \eqref{eq:NewDistances} and present the count of pair distances in Figures \ref{F9}(b),(d),(f),(h). We observe that the analytic and numerical counts are the same in each case. \\

\noindent We next examine the differences between the PCF and the cPCF for a range of domains that contain inaccessible sites. In Figure \ref{F10} we present four domains, containing 0, 1, 25 and 576 clusters of inaccessible sites. The remaining sites are populated with agents (red) at random such that the average occupancy of accessible sites is 20\%. We calculate both the PCF and the cPCF for one hundred identically-prepared realizations and present the functions in Figures \ref{F10}(b),(e),(h),(k) and \ref{F10}(c),(f),(i),(l) for the PCF and the cPCF, respectively. As the accessible sites are populated randomly, there should be no pair correlation present and hence $P(m) = 1$ for all $m$. For the domain with 0 inaccessible sites (Figure \ref{F10}(a)), we expect that the PCF and cPCF will be identical, as the cPCF reduces to the PCF in this case. As expected, we observe that the correlation functions in Figures \ref{F10}(b) and \ref{F10}(c) are indistinguishable. For all three domains with inaccessible sites, the PCF incorrectly suggests that pair correlation is present within the images for a range of pair distances. For the domain in Figure \ref{F10}(d), the PCF implies that there is a mechanism that results in both short and long-range clustering, as the correlation value is greater than one for $m < 10$ and $m > 40$. Further, there appears to be a mechanism which inhibits clustering at intermediary distances. A similar trend is observed for the third domain (Figure \ref{F10}(g)), as well as regular oscillations in correlation for $m < 40$. For the fourth domain (Figure \ref{F10}(j)), these oscillations dominate the PCF. In contrast, the cPCF is approximately one for all pair distances and correctly indicates that there is no mechanism influencing clustering present in the locations of the agents. Hence, the correlation observed in Figures \ref{F10}(e),(h),(k) is an artifact associated with the inaccessible sites.\\

\begin{figure}
\begin{center}
\includegraphics[width=1.0\textwidth]{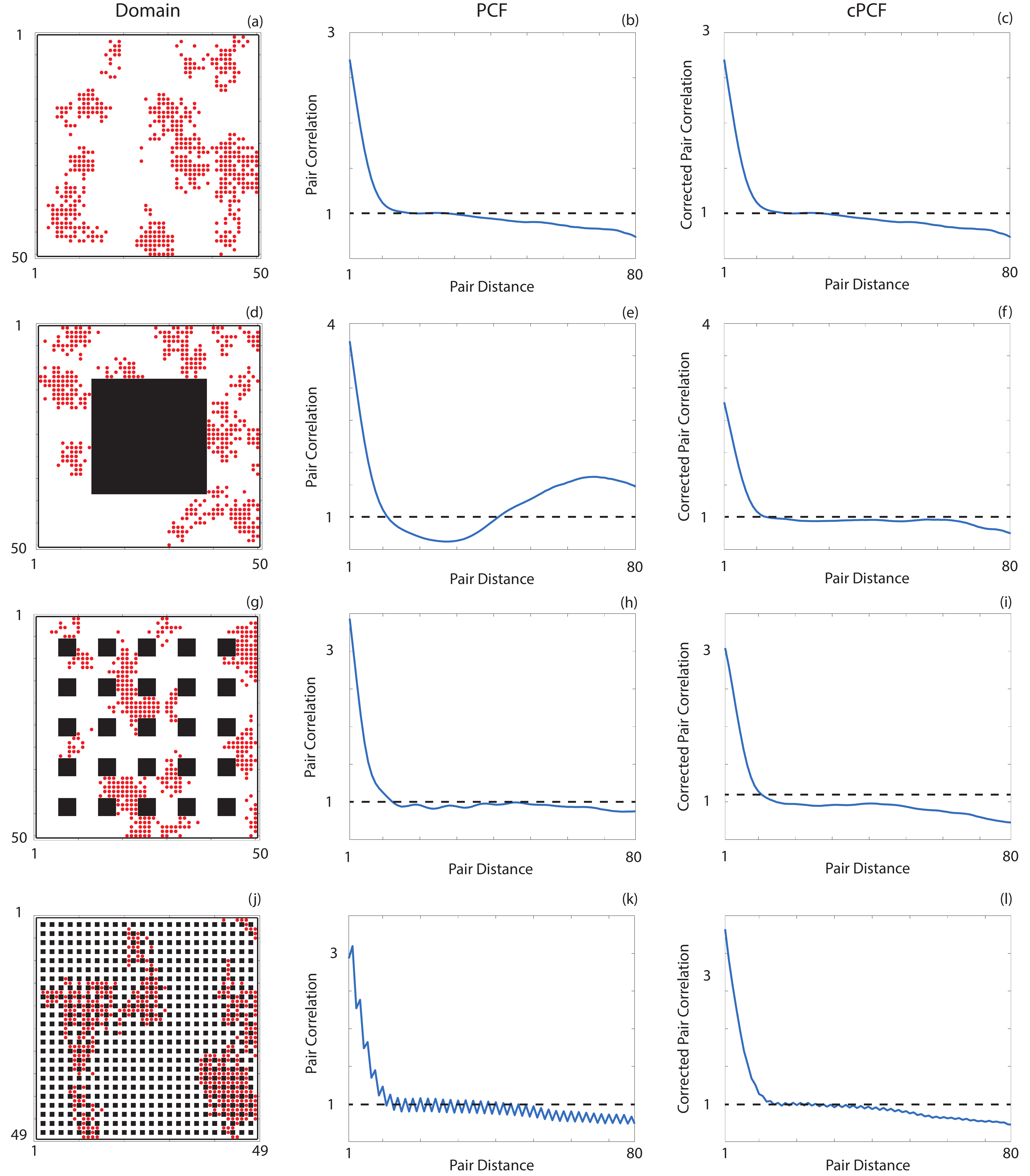}
\caption{(a),(d),(g),(j) Domains containing various configurations of inaccessible sites (black) with agents (red) located at accessible sites after undergoing a birth-movement random walk with the corresponding (b),(e),(h),(k) standard PCF, that is, the pair correlation calculated ignoring inaccessible sites, and (c),(f),(i),(l) cPCF. The dashed black line corresponds to no correlation. All PCFs are the average of 100 identically-prepared domains and the subsequent realization of the random walk process.}
\label{F11}
\end{center}
\end{figure}

\noindent We next compare the PCF and the cPCF for the same domains as analyzed previously, but for agents that follow a birth-movement exclusion-based random walk process. Such processes are discussed in detail elsewhere and have been used to mimic the behavior of a cell population \cite{johnston2014,simpson2010}. Briefly, we populate accessible sites with $z$ agents at random such that the average occupancy is 1\%. Agents undergo birth and movement events with probabilities $P_b$ and $P_m$ per time step, respectively. During a time step, $z$ agents are selected randomly with replacement and undergo birth events, where a daughter agent is placed at one of four randomly selected nearest-neighbor sites \cite{simpson2010}. This birth event is successful if the selected nearest-neighbor site does not contain an agent, and the selected site is not inaccessible. After the birth events have been attempted, $z$ agents are selected randomly with replacement to undergo movement events. During a movement event, one of four nearest-neighbor sites is selected at random, and an agent attempts to move to that lattice site \cite{simpson2010}. Similar to the birth events, this event is successful if the selected site does not contain an agent and the selected site is not inaccessible. This random walk process is used as it is known to result in short-range correlation between agents due to the birth mechanism \cite{johnston2014}, as birth events inherently cause agents to be located at neighboring sites. This clustering tendency is countered by the movement mechanism, which acts in a diffusive manner. Hence, for higher ratios of $P_b$ to $P_m$, we expect to see more short-range correlation, and for $P_m \gg P_b$, we expect to see no correlation. \\

\noindent Four representative snapshots of domain occupancy for agents following the birth-movement random walk process are shown in Figures \ref{F11}(a),(d),(g),(j). In each simulation we use a final time that is weighted by the chance of successfully undergoing movement or birth, as the location and number of the inaccessible sites can influence this chance, and hence comparisons between simulations may have an effectively different final time. For the domain in Figure \ref{F11}(j), for example, there are many sites that have only two accessible neighbor sites. This can be compensated for by scaling the final time by the ratio of the number of accessible neighbor sites if all neighbor sites are accessible to the actual number of accessible neighbor sites in the domain. Hence, for Figure \ref{F11}(j), we scale the final time by approximately 1.48, as there are 5196 out of a possible 7696 accessible neighbor sites. For all simulations $P_m = P_b = 0.1$ and $t_{\text{end}} = 70$, before scaling. Compared to the randomly occupied domains presented in Figures \ref{F10}(a),(d),(g),(j) we immediately observe that the agents are located in clusters. When we calculate the PCFs for these domains, we therefore expect to observe pair correlation values greater than one for short distances. We present the PCF for all four domains in Figures \ref{F11}(b),(e),(h),(k). While we observe that the pair correlation is greater than one at short distances in all four cases, the correlation at longer pair distances varies. For the domain with no inaccessible sites (Figure \ref{F11}(a)), the pair correlation is below one for $m > 20$ and above one for $m < 20$, as expected. In the second domain, the correlation is below one for intermediate pair distances and above one for $m > 40$. The correlation for the third domain is approximately one for intermediary and large $m$ values. For the fourth domain, the pair correlation is only above one for short distances, and below one otherwise. However, there is a distinct oscillatory pattern between odd and even pair distances. As the mechanisms in the random walk process only result in correlation at a pair distance of one, it is unlikely that these oscillations of this magnitude arise from the random walk process. To examine whether these correlations are indeed present due to the random walk mechanisms, we present the cPCF in Figures \ref{F11}(c),(f),(i),(l) for the four domains. Again, for the domain with no inaccessible sites, the PCF and the cPCF are the same. In all cases, we observe the expected high correlation at short distance associated with birth events. For the domain in Figure \ref{F11}(d) the correlation is approximately one for the remainder of the pair distances. For the domains in Figures \ref{F11},(g),(j), however, it appears that the correlation decreases with pair distance. This suggests that the restricted geometry of the domain may influence the spreading of the agents, even while scaling the final time. Interestingly, this decrease is also present in the PCF, for the domain in Figure \ref{F11}(j), albeit in the presence of oscillations. For all four domains the cPCF is relatively consistent, compared to the PCF, which is strongly domain dependent. As such, the cPCF provides a meaningful measure of the correlation present in the domain, as it is able to isolate the agent-agent correlation from spurious correlations arising from domain geometry. \\

\begin{table}
\begin{center}
\begin{tabular}{|r|c|c|c|c|c|c|}
\hline Domain & 50x50 & 50x50 & 100x100 & 100x100 & 150x150 & 150x150 \\ \hline
Number of clusters & 25 & 100 & 25 & 100 & 25 & 100 \\ \hline
Analytic time (s) & 2.18 & 2.36 & 55.61 & 30.35 & 211.64 & 194.83 \\ \hline
Numerical time (s) & 6.48 & 6.41 & 310.85 & 75.58 & 738.05 & 573.98 \\ \hline
\end{tabular}
\caption{Average time taken to evaluate the cPCF for a randomly occupied domain at 20\% density with the specified number of clusters of inaccessible sites for the analytic expression and the numerical path-finding algorithm. Times are reported are the average time taken for 100 randomly generated domains.}
\label{T1}
\end{center}
\end{table}

\begin{table}
\begin{center}
\begin{tabular}{|r|c|c|c|c|c|c|}
\hline Domain & 50x50 & 50x50 & 100x100 & 100x100 & 150x150 & 150x150 \\ \hline
Number of clusters & 25 & 100 & 25 & 100 & 25 & 100 \\ \hline
Analytic time (s) & 3.81 & 3.16 & 62.24 & 56.23 & 345.20 & 352.57 \\ \hline
Numerical time (s) & 19.02 & 14.41 & 348.06 & 304.75 & 2053.61 & 2013.19 \\ \hline
\end{tabular}
\caption{Average time taken to evaluate the cPCF for a randomly occupied domain at 20\% density with the specified number of clusters of inaccessible sites for the analytic expression and the numerical path-finding algorithm. Here the maximum number of inaccessible sites is limited to 25\% of the total number of sites. Times are reported are the average time taken for 100 randomly generated domains.}
\label{T2}
\end{center}
\end{table}

\noindent Finally, we compare the time taken to evaluate the cPCF both using the path-finding algorithm described previously and the analytic expression \eqref{eq:NewDistances}. We consider domains of different sizes randomly occupied by agents such that 20\% of the accessible sites contain an agent for different numbers of clusters of inaccessible sites. In Table \ref{T1} we present the time required to evaluate the cPCF if the domain is randomly generated with no restriction on the number of inaccessible sites, and in Table \ref{T2} we present the time required if the maximum number of inaccessible sites is restricted to 25\% of the total number of sites. In both cases, we observe that the numerical approach is always slower, and requires between three and six times the computation time of the analytic approach. Interestingly, the numerical approach requires additional computation time for fewer inaccessible sites. Furthermore, as the numerical algorithm requires individual calculation of the distance from one site to all other sites for each site (approximately $(L_xL_y)^2$ algorithm realizations), memory issues become a considerable problem as the domain increases in size. 

\subsection*{Approximation}

\begin{figure}
\begin{center}
\includegraphics[width=1.0\textwidth]{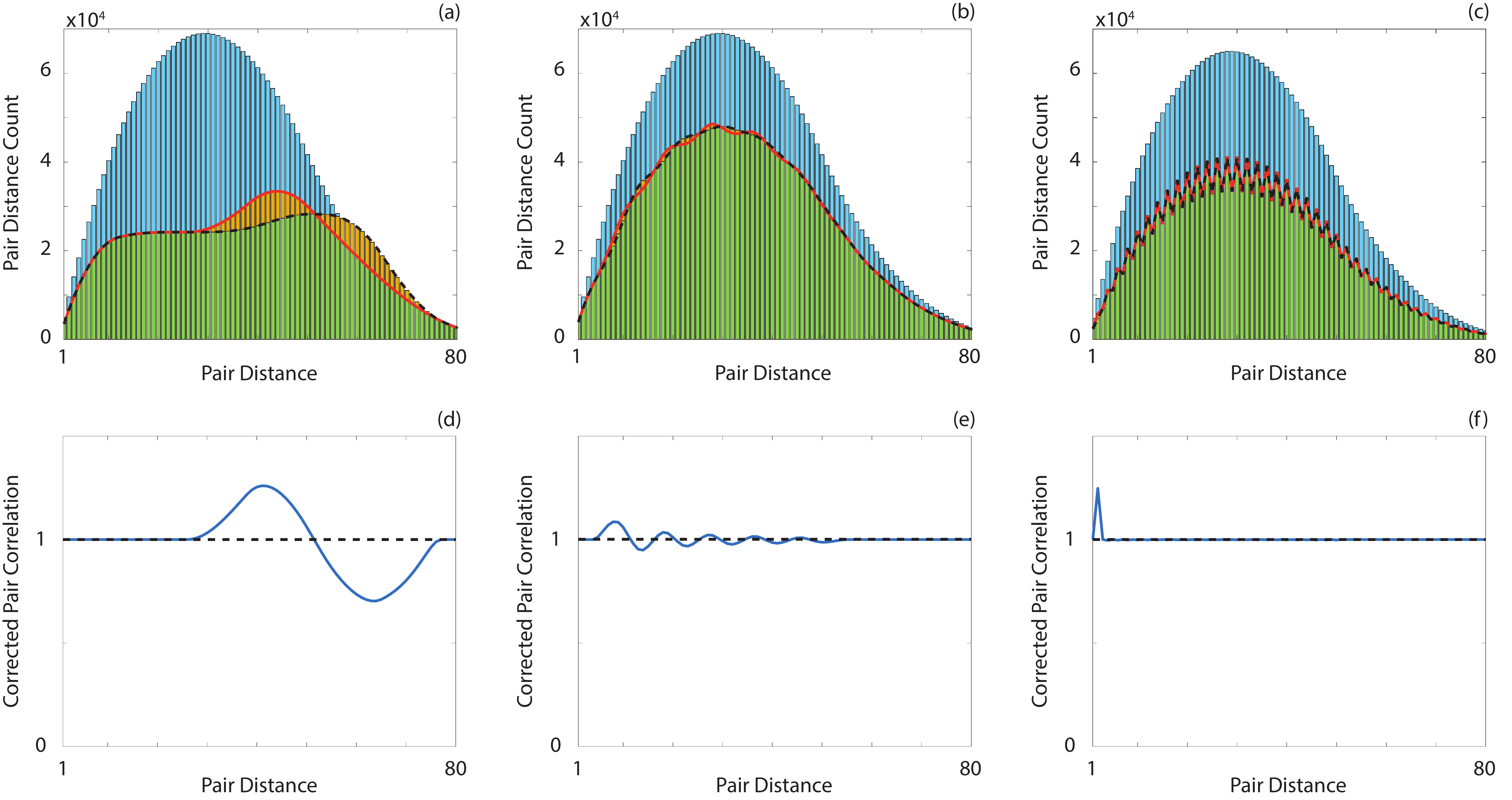}
\caption{(a)-(c) Pair distance count corresponding to the domains presented in Figure \ref{F11}(a),(d),(g), respectively. Blue corresponds to pair distances obtained solely from $D_{\text{NO}}(m)$ (i.e., uncorrected), green corresponds to pair distances from $D_{\text{NO}}(m)$ corrected by $I(m)$ and $A(m)$, and orange represents the correction associated with $S(m)$. The corrected count of pair distances and the approximation are superimposed by the dashed black and red lines. (d)-(f) The approximate cPCF (blue) and exact cPCF (dashed black) for randomly located agents on the domains presented in Figure \ref{F11}(a),(d),(g), respectively.}
\label{F12}
\end{center}
\end{figure}

The cPCF relies on a normalization term, \eqref{eq:NewDistances}, that is only exact for certain configurations of inaccessible sites. As it is computationally intensive to determine the normalization term through many realizations of a path-finding algorithm \cite{gavagnin2018}, which is required for the exact counts of pair distances for such configurations, it is of interest to examine whether an approximation of the analytic normalization term provides a sufficiently accurate alternative. As discussed previously, the corrected normalization term is composed of the standard pairs, accessible-inaccessible pairs, inaccessible-inaccessible pairs and shifted pairs terms. The restriction to certain configurations of inaccessible sites is solely due to the shifted pairs, as the number of other pairs are calculated using standard distance metrics rather than path distance. Hence, if we do not consider the shifted pair distances, the restriction on inaccessible site configurations can be relaxed. We first examine the contributions of shifted pairs to the overall pair distances to determine the size of this contribution. As the shifted pairs consist of both negative and positive terms, corresponding to lost and gained pairs, respectively, the combined terms may only provide a small contribution to the total pairs. We note that for each lost pair there is a corresponding gained pair and hence the total number of pairs is constant independent of whether the lost and gained pairs are considered. In Figure \ref{F12} we present the contribution of the shifted pairs term to the overall count of pair distances for the three domains in Figure \ref{F11}(d),(g),(j). The blue bars correspond to the standard pair distance counts, $D_{\text{NO}}(m)$, the green bars correspond to standard pair distance counts corrected by accessible-inaccessible and inaccessible-inaccessible pairs, and the orange bars correspond to the correction associated with the shifted pairs. As such, the red line corresponds to the approximation of the pair distance count and the black dashed line corresponds to the exact corrected pair distance count. We observe that the shifted pairs provide a small contribution, except in the case of a single large cluster of inaccessible sites (Figure \ref{F12}(a)). As such, an approximation of the corrected normalization term may result in a valid approximation of the cPCF, provided that the domain is not dominated by a single large cluster of inaccessible sites. We note that for these domains we are able to compare the analytic pair distances with the approximation as the configuration of inaccessible sites means that the analytic function is exact. The approximation of the counts of the pair distances is
\begin{equation}
D_{\text{approx}}(m) = D_{\text{NO}}(m) - A(m) + I(m),
\label{eq:ApproxDistances}
\end{equation}
where $D_{\text{NO}}(m)$, $A(m)$ and $I(m)$ are defined in \eqref{eq:Distances}, \eqref{eq:A} and \eqref{eq:I}, respectively.  \\

\begin{figure}
\begin{center}
\includegraphics[width=1.0\textwidth]{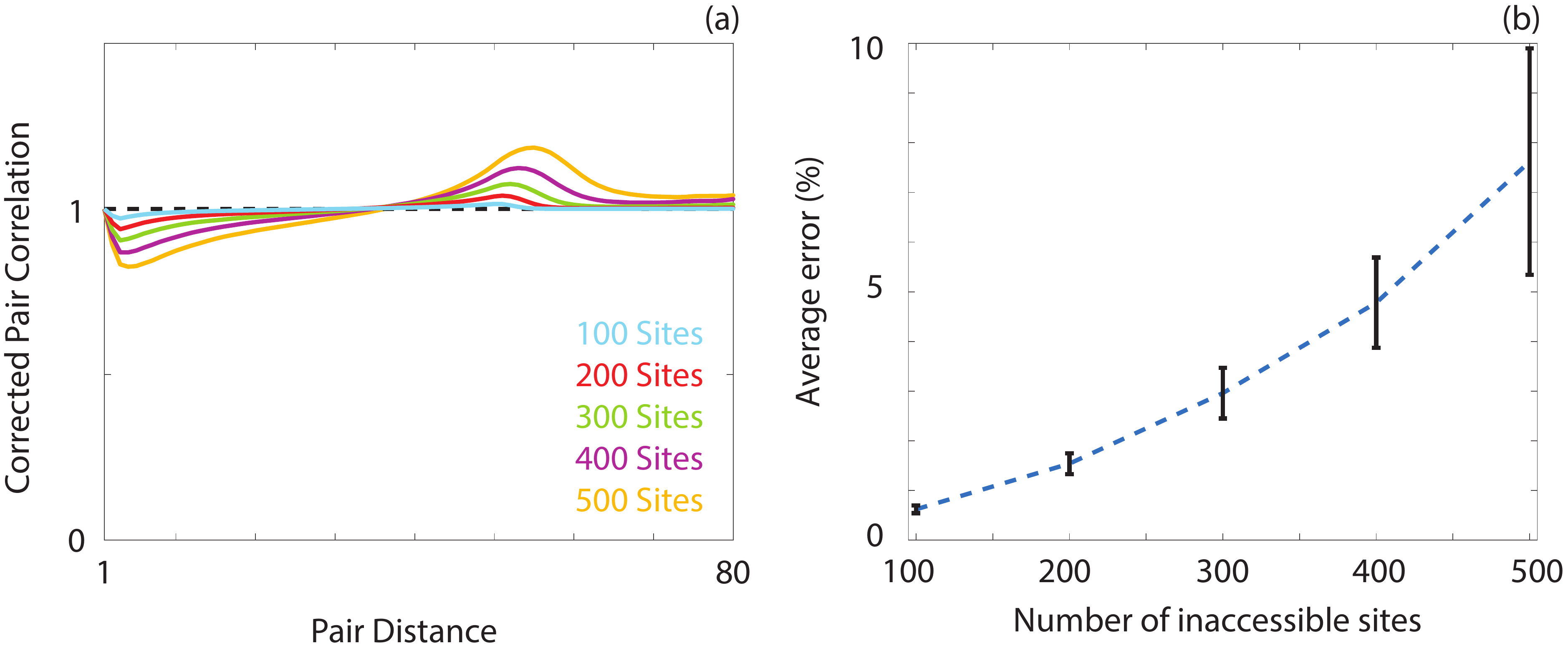}
\caption{(a) Approximate cPCF for 50 randomly generated domains containing 100 (cyan), 200 (red), 300 (green), 400 (purple) or 500 (orange) inaccessible sites for agents placed at random on accessible sites. The exact cPCF is shown via the dashed black line. (b) Error (mean $\pm$ one standard deviation) in the approximate cPCF, compared to the exact result, as a function of the number of inaccessible sites.}
\label{F13}
\end{center}
\end{figure}

\noindent For the domains considered previously we present both the cPCF and the corresponding approximation in Figures \ref{F12}(d)-(f). As expected, we see that for the domain with a large cluster of inaccessible sites, the approximation is poor for pair distances similar in size to the cluster. For the other two domains, the approximate cPCF performs well. Finally, we populate domains with inaccessible sites at random and calculate the approximate cPCF. All accessible sites on the domain are populated by agents such that the expected pair correlation is one for all pair distances. In Figure \ref{F13}, we present the average approximate cPCF for 50 random identically-prepared domains for a range of numbers of inaccessible sites, as well as the mean error associated with each approximate cPCF. Intuitively, we observe that an increase in inaccessible sites corresponds to an increase in the distance between the expected pair correlation and the approximation. Increasing the number of inaccessible sites while populating the domain with inaccessible sites at random introduces accessible sites that are not connected to the remainder of the domain, and hence we do not consider the cPCF approximation for higher numbers of inaccessible sites. As such, the approximation may prove useful for calculating pair correlations for domains where the inaccessible sites do not satisfy the conditions required for the exact pair distance calculation, but have a modest proportion of inaccessible sites compared to accessible sites.

\section{Discussion and conclusions}
\label{s:Discussion}
\noindent Analysis of the spatial structure present in experimental images provides valuable insight into the mechanisms governing behavior within the experiment \cite{johnston2014,law2009,perry2002}. Pair correlation functions have been employed in a variety of fields, including ecology \cite{law2009}, biology \cite{dini2018} and physics \cite{bahcall1983}, and have proven useful for elucidating the presence and impact of spatial structure. Many experimental environments contain immobile obstacles that influence the transport and location of individuals within that environment \cite{rehder2015,roosen2011,moussaid2011}. Isolating the spatial structure associated with the mechanisms that govern transport, rather than the heterogeneous nature of the environment is therefore of interest. Naively applying standard PCFs does not account for distances between pairs of individuals that must avoid obstacles, and may result in the incorrect suggestion of spatial correlations.  \\

\noindent Here we have presented an exact analytic expression for the normalization term of a corrected PCF that incorporates a physical path distance between individuals, and hence can be applied to environments with obstacles. We demonstrate that this cPCF is necessary for isolating the spatial correlation associated with the locations of individuals from the spatial correlation associated with the environment itself. Further, we highlight that the analytic expression allows for the cPCF to be calculated significantly faster than relying on a path-finding algorithm. We apply the cPCF to images arising from a lattice-based movement-birth random walk, which mimics cell motility and cell proliferation, where short-range correlation is known to exist, and demonstrate the cPCF recovers this correlation. Standard PCFs can introduce spurious correlations as well as oscillations in the correlation. Finally, we present an approximation to the cPCF that relaxes assumptions on the locations of the inaccessible sites within the domain and show that for modest numbers of inaccessible sites, the approximation is accurate. \\

\noindent The work and analysis presented here could be extended in a number of directions. One obvious application is to calculate the PCF for experimental data obtained from an environment that contains obstacles. Here we have focused on data arriving from simulations that mimic processes such as cell migration and proliferation \cite{johnston2014} rather than explicitly using experimental data. As previous investigations involving the application of PCFs to experimental data have proved fruitful \cite{agnew2014,johnston2014,dini2018}, the application of the cPCF to appropriate data may prove to be insightful. Another promising approach would be to examine how the pair correlation changes between two experiments on different domains. As the cPCF is able to isolate the correlation associated with the behavior of individuals, it would be instructive to consider whether the behavior is dependent on environment and, if so, quantify which mechanisms are are responsible for this change in behavior.

\begin{acknowledgments}
This research was in part conducted and funded by the Australian Research Council Centre of Excellence in Convergent Bio-Nano Science and Technology (project no. CE140100036). We thank Pradeep Rajasekhar and Daniel Poole for providing the image of the nervous system used in Figure 1(a).
\end{acknowledgments}


\appendix
\section{Standard counts of pair distances derivation}
\noindent Consider an arbitrary domain with no obstacles and $L_x$ by $L_y$ lattice sites. For pair distances that satisfy $1 \leq m \leq \text{min}(L_x,L_y)$, the number of pairs of sites separated by a distance $m$ was presented by Gavagnin \emph{et al.} \cite{gavagnin2018}. However, the maximum pair distance on such a domain with no-flux boundary conditions is $L_x+L_y-2$, and we therefore must derive counts of pair distances for $m > \text{min}(L_x,L_y)$. \\

\noindent We first consider pair distances where $\text{min}(L_x,L_y) < m < \text{max}(L_x,L_y)$, that is, pairs of sites that cannot be connected via a path containing only ``jumps'' in the shorter of the horizontal or vertical directions. For example, if $L_y < L_x$, the path connecting pairs of sites separated by distance $\text{min}(L_x,L_y) < m < \text{max}(L_x,L_y)$ must contain either entirely horizontal ``jumps'', or a combination of horizontal and vertical ``jumps''. For a particular $m$, the number of pairs of sites separated by $m$ horizontal ``jumps'' is $L_y(L_x-m)$. This term is the product of the number of pairs of sites separated by $m$ horizontal ``jumps'' within a single row, $L_x-m$, and the number of rows, $L_y$. The number of sites separated by $m-1$ horizontal ``jumps'' and 1 vertical ``jump'' is $2(L_y-1)(L_x-m+1)$. Compared to the $m$ horizontal ``jumps'' case, there is an additional number of pairs of sites separated by $m-1$ horizontal ``jumps'', providing the $L_x-m+1$ component. However, as each distance contains a vertical ``jump'', the rows containing the sites must be offset by one and hence there is one less row where this can occur, which results in the $L_y-1$ term. Finally, as the vertical ``jump'' can be in either the positive or negative vertical direction, this introduces the factor of two. Introducing additional vertical ``jumps'' increases the number of pairs of sites separated by the resulting smaller number of horizontal ``jumps'', while decreasing the effective number of rows. Noting that $L_y$ and $L_x$ are interchangeable, we therefore obtain
\begin{align*}
D_{\text{NO}}(m) & = \text{min}(L_x,L_y)\bigg[\text{max}(L_x,L_y)-m\bigg] \\ & + 2\sum_{i=1}^{\text{min}(L_x,L_y)}\bigg[\text{min}(L_x,L_y)-i\bigg]\bigg[\text{max}(L_x,L_y)-m+i\bigg], \\ & = D_{\text{NO}}(\text{min}(L_x,L_y)) - \text{min}(L_x,L_y)^2(m-\text{min}(L_x,L_y))
\end{align*}
for $\text{min}(L_x,L_y) < m < \text{max}(L_x,L_y)$. \\

\noindent We next consider pair distances $m \geq \text{max}(L_x,L_y)$, where distances between sites must include both horizontal and vertical ``jumps''. Again, without loss of generality, we assume that $L_y \leq L_x$. The minimum number of vertical ``jumps'' for a distance $m$ is $m-L_x+1$, and the maximum number is $L_y$. The corresponding number of pairs of rows that are separated by this vertical distance $v$ is $L_y-v$, and the number of pairs of sites separated by this vertical distance and the requisite horizontal distance $h = m-v$ is $L_x-m+v$. Following the same process as above, taking a summation over the possible rows and columns, and noting that the vertical separation can be either in the positive or negative direction, we obtain
\begin{align*}
D_{\text{NO}}(m) & = 2\sum_{i=m-\text{max}(L_x,L_y)+1}^{\text{min}(L_x,L_y)} \bigg[\text{min}(L_x,L_y)-i\bigg]\bigg[\text{max}(L_x,L_y)-m+i\bigg], \\ & = \frac{k(k+1)(k+2)}{3}, \ \text{where} \ k = L_x+L_y-1-m,
\end{align*}
for $m \geq \text{max}(L_x,L_y)$. 

\bibliographystyle{apsrev}

\end{document}